\documentclass{article}
\usepackage{bbm}
\usepackage{amssymb}
\usepackage{amsmath}
\usepackage{amsfonts}
\usepackage{mathrsfs}
\usepackage {amssymb}
\usepackage {amsmath}
\usepackage {amsthm}
\usepackage{graphicx,color}
\usepackage {amscd}
\usepackage{geometry}
\usepackage{hyperref}

\title{Regularity for Harmonic - Einstein Equation}
\date{\today}
\author{Yiyan Xu \footnote{\href{mailto:xuyiyan@math.pku.edu.cn}{\textcolor[rgb]{0.00,0.00,1.00}{Email: xuyiyan@math.pku.edu.cn}}}}
\theoremstyle{definition}
\newtheorem{lemma}{Lemma}[section]
\newtheorem{definition}[lemma]{Definition}

\newtheorem{theorem}[lemma]{Theorem}

\newtheorem{remark}{Remark}

\numberwithin{equation}{section}
\renewcommand{\proof}{Proof. }
\renewcommand{\qed}{\hfill\small{$\square$}\normalsize}
\begin{document}
\maketitle
\begin{abstract}
We establish a regularity theorem for the Harmonic - Einstein
Equation. As a byproduct of the local regularity, we also have a
compactness theorem on Harmonic - Einstein equation.
The method is mainly the Moser iteration technique which has been
used and developed by \cite{BKN89}, \cite{Tian90}, \cite{TV05a} and
others.
\end{abstract}
\section{Introduction}
In this paper, we will consider the degeneration of Harmonic -
Einstein equation,  which is a generalization the Einstein equation.
Another motivation is that John Lott  showed that the expanding
soliton equation on the space with the simplest type of Nil
structure can be reduced to the Harmonic - Einstein equation, while
such kind of soliton appeared in the long time limit of type - III
Ricci flow.

\theorem \cite{Lott07} Let $(\overline{M},\overline{g})$ be the
total space of a flat $\mathbb{R}^N$- vector bundle over a
Riemannian manifold $(M,g)$, with flat Riemannian metrics on the
fibers. Suppose that the fiberwise volume forms are preserved by the
flat connection. Let $V$ be the fiberwise radial vector field
$\frac{1}{2}\sum_{i=1}^Nx^i\frac{\partial}{\partial x_i}$. Then the
expanding soliton equation on $\overline{M}$
\[\overline{Ric}+\frac{1}{2}\mathscr{L}_V\overline{g}+\frac{1}{2}\overline{g}=0\]
becomes the equation for a
harmonic map
\begin{equation}\label{HaME}
    G: (M,g)\rightarrow (SL(N,\mathbb{R})/SO(N),\langle\;,\;\rangle)
\end{equation}
along the equation
\begin{equation}\label{HEM}
   Ric-\frac{1}{4}\langle dG,dG\rangle+\frac{1}{2}g=0
\end{equation}
on $M$, where $\langle\;,\;\rangle_G=Tr(G^{-1}dGG^{-1}dG)$ is the
usual metric on the symmetric space $SL(N,\mathbb{R})/SO(N)$.
\remark There is an algebraic description of Symmetric space
$SL(N,\mathbb{R})/SO(N)$: $sl(N)=\{X:trX=0\}\simeq \hbar\oplus
so(N)$, where $\hbar$ is the symmetric part, which can be identified
with the tangent space of $SL(N,\mathbb{R})/SO(N)$. On $\hbar$ we
have the usual Euclidean metric, the involution $L$ is $-id$ on
$\hbar$ and $id$ on $so(N)$, namely, $L(X)=-X^t.$ Consequently, the
curvature is $Rm(X,Y,Z, W)=-\langle[X,Y],[Z,W]\rangle.$  In
particular, the sectional curvature is nonpositive, which is crucial
in our result.
\definition Let $G:(M,g)\longrightarrow(N,h)$ be a map between two Riemannian
manifolds, with the notation  in \cite{Lott07},  we will call the
local version of the equations
\begin{equation}\label{HaMEl}
\left\{ \begin{aligned}
         0&= Ric(g)-\langle dG, dG\rangle-\lambda g \\
          0&=\Delta_{g,h} G
         \end{aligned} \right.
\end{equation}
or in local coordinates,
\begin{equation}
\left\{ \begin{aligned}
         0&= Ric(g)_{\alpha\beta}-h_{ij}G_{,\alpha}^iG_{,\beta}^j-\lambda g_{\alpha\beta} \\
          0&=g^{\alpha\beta}G_{,\alpha\beta}^i-g^{\alpha\beta}G_{,\alpha}^jG_{,\beta}^k\Gamma(h)_{jk}^i
         \end{aligned} \right.
\end{equation}
 to be \texttt{Harmonic-Einstein} equation.

Since any Riemannian metric satisfies $\Delta Rm=\nabla^2Ric+Rm\ast
Rm$, combine with the equation (\ref{HaMEl}), $(M,g,G)$ satisfies a
coupled Elliptic system, together with uniform Sobolev constant
$C_S$, one can prove an $\epsilon$ - regularity theorem.

\theorem\label{eRegT} ($\epsilon$ - regularity) Assume $(M,g,G)$
satisfies the Harmonic-Einstein equation (\ref{HaMEl}), and $(N,h)$
has nonpositive sectional curvature, $\lambda=0$ or $-1$. Let
$B(0,r)$ be a geodesic ball around $0\in M$, $C_S$ be the Sobolev
constant on $B(0,r)$, and $k\in\mathbb{N}$. Then there exists a
constant $\epsilon=\epsilon(C_S,n)$ such that if
\[\{\int_{B(0,r)} |Rm|^{\frac{n}{2}}\}^{\frac{2}{n}}\leq\epsilon,\]
then
\begin{eqnarray*}
 \sup_{B(0,\frac{r}{2})}|\nabla^{k}
 G|&\leq&\frac{C}{r^{k}}(\{\int_{B(0,r)}|Rm|^{\frac{n}{2}}\}^{\frac{2}{n}}+\{r^{2-n}\int_{B(0,r)}|d
  G|^2\}^{\frac{1}{2}}),\\
  \sup_{B(0,\frac{r}{2})}|\nabla^kRm|&\leq&
  \frac{C}{r^{k+2}}(\{\int_{B(0,r)}|Rm|^{\frac{n}{2}}\}^{\frac{2}{n}}+\{r^{2-n}\int_{B(0,r)}|d
  G|^2\}^{\frac{1}{2}}),
\end{eqnarray*}
where $C=C(C_S,k,n)$.

At this point, the author would like to point out the method here
can also be used to prove the $\epsilon$ - regularity theorem on the
system $\Delta Ric= Rm\ast Ric$. The four dimensional case, namely,
Bach fat metric with constant scalar curvature, has been established
by Tian and Viaclovsky \cite{TV05a}, and the higher dimension case
has been proved by Chen and Weber \cite{CW11}. The main idea are all
similar, but one will see our iteration process is different.

As a byproduct of the $\epsilon$ - regularity, we obtain a
convergence theorem for Harmonic - Einstein equation, which is
similar to the compactness on harmonic maps \cite{SaUh81}, Yang -
Mills connections \cite{Uhlenbeck82a}, Einstein metrics
\cite{Anderson89}, \cite{BKN89}, \cite{Tian90}, \cite{Nakajima94},
and more recently Bach flat metric with constant scalar curvature
\cite{TV05b}, \cite{AAJV11}, K\"{a}hler Ricci soliton \cite{CS07},
extremal K\"{a}hler metric \cite{CW11}.

\theorem\label{Compactness}(Compactness) Let $(g_i, G_i)$  satisfy
the Harmonic - Einstein equation (\ref{HaMEl}) over a sequence of
$4$ - dimensional compact manifolds $M_i$, respectively. Assume
$(M_i,g_i)$ satisfy:
\[\hbox{Euler number} \; \mathcal {X}(M_i)\leq X,\; Diam (M_i, g_i)\leq D,\; Vol (M_i, g_i) \geq
V,\] and $G_i: (M_i,g_i)\longrightarrow (N,h)$ with finite energy
\[E(G_i, g_i):= \int_{M_i}|dG_i|^2\leq E,\;\]
where $X$, $E$, $D$, $V$ are constants which are independent of $i$.
We also assume $(N,h)$ has nonpositive sectional curvature and
$\lambda=0$ or $-1$. Then there exists a subsequence
$\{j\}\subset\{i\}$ satisfies the following properties:
\begin{enumerate}
  \item $\{M_j, g_j,G_j\}$ converges to a complete metric space $M_\infty$ in the following sense:
  If we remove a finite set $\mathcal {S}=\{b_1,\cdots,b_m\}\subset
  M_\infty$ with $m\leq m(n,X,D,V)$, a $C^\infty$ manifold structure is defined and also a smooth pair $(g_\infty,G_\infty)$
  satisfies the
  Harmonic -
  Einstein equation over the punctured set $M_\infty\setminus\mathcal
  {S}$. Moreover, there exists a (into) diffeomorphism $F_j: M_\infty\setminus\mathcal
  {S}\hookrightarrow M_i$ such that $(F_j^*g_j, F_j^*G_j)$
  converges to $(g_\infty,G_\infty)$ in the $C^\infty(M_\infty\setminus\mathcal
  {S})$ topology.
  \item The manifold structure and the pair $(g_\infty,G_\infty)$ on $M_\infty\setminus\mathcal
  {S}$ extend to the whole of $M$ which satisfies the
Harmonic - Einstein equation (\ref{HaMEl}) over a Riemannian
orbifold.
\end{enumerate}
\definition By $(g,G)$ satisfies the Harmonic -
  Einstein equation (\ref{HaMEl}) over a Riemannian orbifold, $(M,g)$, we mean:
  \begin{enumerate}
    \item There exists a finite set $\mathcal {S}=\{b_1,\cdots,b_m\}\subset
  M$, such that $M\setminus\mathcal{S}$ is a $C^\infty$ manifold and the
    restriction of $(g,G)$ satisfies the smooth Harmonic -
  Einstein equation (\ref{HaMEl}).
    \item For each singular point $b_k\in \mathcal {S}$, there exists a neighborhood $U_k\subset M$ such that $U_k\setminus\{x_k\}$ is
  diffeomorphic to $B^n\setminus\{0\}/\Gamma_k$, where $B^n\in \mathbb{R}^n$ is
  a $n$ - dimensional unit ball and $\Gamma_k\subset O(n)$ is a finite
  subgroup acting freely on $B^n\setminus\{0\}/\Gamma$. If we lift $(g,G)$
  to $B^n\setminus\{0\}$, it extends smoothly across the singular point
  0 and satisfies the Harmonic -
  Einstein equation (\ref{HaMEl}) over $B^n$.
 \end{enumerate}
\textbf{Acknowledgements}: The author would like to thank his
advisor Gang Tian for suggesting this problem.
\section{Bochner Identities}
Let $f:(M,g)\longrightarrow(N,h)$ be a map between two Riemannian
manifolds, the differential of $f$ is
\[df=\frac{\partial f^i}{\partial x^\alpha}dx^\alpha\otimes\frac{\partial}{\partial f^i}\in \Gamma(T^*M\otimes f^{-1}TN).\]
Now we would like to establish the Bochner - type identity on the
extended bundle $\Omega^pT^*M\otimes\Omega^qf^{-1}TN$ over $M$ with
respect to the induced metric $g\otimes f^*h$ and induced connection
$\nabla^M\otimes f^*\nabla^N$. These kind identities should be well
known by experts.

\lemma\label{BTIF} Let $f:(M,g)\longrightarrow(N,h)$ be a map
between two Riemannian manifolds, then we have the well known
Bochner formula, for instance, see \cite{EL78} and \cite{SY94}.
\[\Delta\frac{1}{2}|df|^2=|\nabla df|^2+\langle\nabla\Delta f,\Delta f\rangle-\sum_{\alpha,\beta}Rm^N(f_*e_\alpha,f_*e_\beta,f_*e_\alpha,f_*e_\beta)+\sum_{i}Ric^M(f^*\theta_i,f^*\theta_i)\]
where $\{e_\alpha\}$ is an orthonormal basis for $TM$,
$\{\theta_i\}$ is an orthonormal basis for $T^*N$. More generally,
we also have the Bochner type identities for commutation of
covariant derivatives up to order $k$:
\begin{eqnarray}
\Delta\nabla^k f&=&\nabla^k\Delta
f+\sum_{i=0}^{k-1}\nabla^iRm^M\ast\nabla^{k-i}f\nonumber\\&&+\sum_{p=3}^{k+2}\sum_{i_1+\cdots+i_p=k-p+2}\nabla^{p-3}Rm^N\ast\nabla^{i_1+1}f\ast\cdots\ast\nabla^{i_p+1}f\label{BIHk}
\end{eqnarray}
In particular, if the target manifold is symmetric, i.e. $\nabla
Rm^N=0$, then we can drop the terms which involve the derivative of
$Rm^N$ in the above expression:
\begin{eqnarray*}
  \Delta\nabla^kf &=&\nabla^k\Delta
f+\sum_{i=0}^{k-1}\nabla^iRm^M\ast\nabla^{k-i}f+\sum_{i=0}^{k-1}\sum_{p+q=i}Rm^N(\nabla^{p+1}f,\nabla^{q+1}f)\nabla^{k-i}f\\
\end{eqnarray*}
\proof In order to simplify the computation, we choose normal
coordinates at $x$ and $f(x)$ respectively, namely,
\[g_{\alpha\beta}(x)=\delta_{\alpha\beta}, g_{\alpha\beta,\gamma}(x)=0;\quad h_{ij}(f(x))=\delta_{ij}, h_{ij,k}(f(x))=0.\]
Therefore, we only have to take the second and up derivatives of the
metric into account, and these will turn into the curvature terms.
First, let us compute the commutation of covariant derivatives up to
three order directly:
\begin{eqnarray*}
 &&\nabla^3 f(\frac{\partial}{\partial x^\alpha},\frac{\partial}{\partial x^\beta},\frac{\partial}{\partial x^\gamma}) =\nabla_{\frac{\partial}{\partial x^\alpha}}(\nabla^2 f(\frac{\partial}{\partial x^\beta},\frac{\partial}{\partial
 x^\gamma}))\\
  &=&\nabla_{\frac{\partial}{\partial x^\alpha}}(\frac{\partial^2f^i}{\partial x^\beta\partial x^\gamma}\frac{\partial}{\partial f^i}+\frac{\partial f^i}{\partial x^\gamma}\frac{\partial f^j}{\partial x^\beta}\nabla_{\frac{\partial}{\partial f^j}}\frac{\partial}{\partial f^i}-df(\nabla_{\frac{\partial}{\partial x^\beta}}\frac{\partial}{\partial x^\gamma}))\\
   &=&\frac{\partial^3f^i}{\partial x^\alpha\partial x^\beta\partial x^\gamma}\frac{\partial}{\partial f^i}+\frac{\partial f^i}{\partial x^\gamma}\frac{\partial f^j}{\partial x^\beta}\frac{\partial f^k}{\partial x^\alpha}\nabla_{\frac{\partial}{\partial f^k}}\nabla_{\frac{\partial}{\partial f^j}}\frac{\partial}{\partial f^i}-df(\nabla_{\frac{\partial}{\partial x^\alpha}}\nabla_{\frac{\partial}{\partial x^\beta}}\frac{\partial}{\partial x^\gamma})\\
   &=&\frac{\partial^3f^i}{\partial x^\alpha\partial x^\beta\partial x^\gamma}\frac{\partial}{\partial f^i}
   +\frac{\partial f^i}{\partial x^\gamma}\frac{\partial f^j}{\partial x^\beta}\frac{\partial f^k}{\partial x^\alpha}(\nabla_{\frac{\partial}{\partial f^i}}\nabla_{\frac{\partial}{\partial f^k}}\frac{\partial}{\partial f^j}-Rm^N(\frac{\partial}{\partial f^k},\frac{\partial}{\partial f^i})\frac{\partial}{\partial f^j})\\&&-df(\nabla_{\frac{\partial}{\partial x^\gamma}}\nabla_{\frac{\partial}{\partial x^\alpha}}\frac{\partial}{\partial x^\beta}-Rm^M(\frac{\partial}{\partial x^\alpha},\frac{\partial}{\partial x^\gamma})\frac{\partial}{\partial x^\beta})\\
   &=&\nabla^3 f(\frac{\partial}{\partial x^\gamma},\frac{\partial}{\partial x^\alpha},\frac{\partial}{\partial x^\beta})-Rm^N(\frac{\partial f^k}{\partial x^\alpha}\frac{\partial}{\partial f^k},\frac{\partial f^i}{\partial x^\gamma}\frac{\partial}{\partial f^i})\frac{\partial f^j}{\partial x^\beta}\frac{\partial}{\partial f^j}+df(Rm^M(\frac{\partial}{\partial x^\alpha},\frac{\partial}{\partial x^\gamma})\frac{\partial}{\partial
   x^\beta})
\end{eqnarray*}
taking trace with respect to $\alpha$ and $\beta$, we obtain
\[\Delta\nabla f=\nabla\Delta f+Rm^N(f_*(\frac{\partial}{\partial x^\alpha}),df)f_*(\frac{\partial}{\partial x^\alpha})+df(Ric^M)\]
Consequently, we get the Bochner formula:
\begin{eqnarray*}
  \Delta\frac{1}{2}|df|^2 &=& |\nabla df|^2+\langle\Delta \nabla f,\nabla f\rangle \\
   &=& |\nabla df|^2+\langle\nabla\Delta f,\nabla f\rangle -Rm^N(f_*(\frac{\partial}{\partial x^\alpha}),f_*(\frac{\partial}{\partial x^\gamma}),f_*(\frac{\partial}{\partial x^\alpha}),f_*(\frac{\partial}{\partial x^\gamma}))\\&&+\langle df(Ric^M(\frac{\partial}{\partial x^\gamma})),df(\frac{\partial}{\partial x^\gamma})\rangle
\end{eqnarray*}
Therefore, we proved the case $k=1$. Now let us assume the
expression (\ref{BIHk}) holds for $k-1$, and we will prove the case
$k$.

More generally, on the extended bundle
$\Omega^pT^*M\otimes\Omega^qf^{-1}TN$ over $M$ with respect to the
induced connection $\nabla^M\otimes f^*\nabla^N$, for
 \[T=T_{\alpha_1,\cdots,\alpha_p}^{i_1,\cdots,i_q}dx^{\alpha_1}\otimes\cdots\otimes
dx^{\alpha_p}\otimes\frac{\partial}{\partial
f^1}\otimes\cdots\otimes\frac{\partial}{\partial f^q}\in
\Omega^pT^*M\otimes\Omega^qf^{-1}TN,\] with abuse of notation:
\[\nabla^2_{\frac{\partial}{\partial x^\alpha},\frac{\partial}{\partial x^\beta}}T-\nabla_{\frac{\partial}{\partial x^\beta},\frac{\partial}{\partial x^\alpha}}^2T=R(\frac{\partial}{\partial x^\alpha},\frac{\partial}{\partial x^\beta})T\]
where
\begin{eqnarray*}
Rm(T) &=&
\sum_{i=1}^pT_{\alpha_1,\cdots,\alpha_p}^{i_1,\cdots,i_q}dx^{\alpha_1}\otimes\cdots\otimes
Rm^M(dx^{\alpha_i})\otimes\cdots\otimes
dx^{\alpha_p}\otimes\frac{\partial}{\partial
f^1}\otimes\cdots\otimes\frac{\partial}{\partial
f^q} \\
&&
+\sum_{j=1}^qT_{\alpha_1,\cdots,\alpha_p}^{i_1,\cdots,i_q}dx^{\alpha_1}\otimes\cdots\otimes
dx^{\alpha_p}\otimes\frac{\partial}{\partial f^1}\otimes \cdots
\otimes Rm^N(df, df)(\frac{\partial}{\partial f^j})\otimes\cdots
\otimes\frac{\partial}{\partial f^q}\\
&=&(Rm^M+Rm^N(df,df))\ast T
\end{eqnarray*}
then
\begin{eqnarray}
\nabla Rm(T) &=&\nabla\Big((Rm^M+Rm^N(df,df))\ast T\Big)\nonumber\\
&=&(\nabla Rm^M+Rm^N(\nabla df,df)+df\otimes\nabla Rm^N(df,df))\ast
T\nonumber\\&&+(Rm^M+Rm^N(df,df))\ast \nabla T\label{CCDFR}
\end{eqnarray}
With this notation, let us compute the $\Delta\nabla^kf$,
\begin{eqnarray*}
  && \Delta\nabla^kf(\frac{\partial}{\partial x^{\beta_1}},\cdots,\frac{\partial}{\partial x^{\beta_k}})  \\
   &=& \nabla^{k+2}f(\frac{\partial}{\partial x^{\alpha}},\frac{\partial}{\partial x^{\alpha}},\frac{\partial}{\partial x^{\beta_1}},\cdots,\frac{\partial}{\partial x^{\beta_k}}) \\
&=&\nabla_{\frac{\partial}{\partial x^{\alpha}}}\Big(\nabla^{k+1}f(\frac{\partial}{\partial x^{\beta_1}},\frac{\partial}{\partial x^{\alpha}},\cdots)+Rm(\nabla^{k-1}f)(\frac{\partial}{\partial x^{\alpha}},\frac{\partial}{\partial x^{\beta_1}},\cdots)\Big)\\
&=&\nabla^{k+2}f(\frac{\partial}{\partial x^{\alpha}},\frac{\partial}{\partial x^{\beta_1}},\frac{\partial}{\partial x^{\alpha}},\cdots)+\nabla(Rm(\nabla^{k-1}f))(\frac{\partial}{\partial x^{\alpha}},\frac{\partial}{\partial x^{\alpha}},\frac{\partial}{\partial x^{\beta_1}},\cdots)\\
&=&\nabla^{k+2}f(\frac{\partial}{\partial x^{\beta_1}},\frac{\partial}{\partial x^{\alpha}},\frac{\partial}{\partial x^{\alpha}},\cdots)+Rm(\nabla^{k}f)(\frac{\partial}{\partial x^{\alpha}},\frac{\partial}{\partial x^{\beta_1}},\frac{\partial}{\partial x^{\alpha}},\cdots)\\&&+\nabla(Rm(\nabla^{k-1}f))(\frac{\partial}{\partial x^{\alpha}},\frac{\partial}{\partial x^{\alpha}},\frac{\partial}{\partial x^{\beta_1}},\cdots)\\
\end{eqnarray*}
From the above calculation, 
by induction and formula (\ref{CCDFR}),
\begin{eqnarray*}
\Delta\nabla^kf&=&\nabla\Delta\nabla^{k-1}f+Rm(\nabla^{k}f)+\nabla(Rm(\nabla^{k-1}f)) \\
&=&\nabla\nabla^{k-1}\Delta
f+\sum_{i=0}^{k-2}\nabla\Big(\nabla^iRm^M\ast\nabla^{k-i-1}f\Big)+(Rm^M+Rm^N(df,df))\ast\nabla^{k}f\\&&+\sum_{p=3}^{k+1}\sum_{i_1+\cdots+i_p=k-p+1}\nabla\Big(\nabla^{p-3}Rm^N\ast\nabla^{i_1+1}f\ast\cdots\ast\nabla^{i_p+1}f\Big)\\
&&+(\nabla Rm^M+Rm^N(\nabla df,df)+df\otimes\nabla
Rm^N(df,df))\ast\nabla^{k-1}f\\
&=&\nabla^k\Delta
f+\sum_{i=0}^{k-1}\nabla^iRm^M\ast\nabla^{k-i}f\\&&+\sum_{p=3}^{k+2}\sum_{i_1+\cdots+i_p=k-p+2}\nabla^{p-3}Rm^N\ast\nabla^{i_1+1}f\ast\cdots\ast\nabla^{i_p+1}f
\end{eqnarray*}
Thus the formula (\ref{BIHk}) holds for $k$.

In particular, if the target manifold is a symmetric space, i.e.
$\nabla Rm^N=0$, therefore we can drop the terms which involve the
derivative of $Rm^N$ in the above expression. With this
simplification, we have the Bochner identity:
\begin{eqnarray*}
  \Delta\nabla^kf &=&\nabla^k\Delta
f+\sum_{i=0}^{k-1}\nabla^iRm^M\ast\nabla^{k-i}f+\sum_{i=0}^{k-1}\sum_{i_1+i_2=i}Rm^N\ast\nabla^{i_1+1}f\ast\nabla^{i_2+1}f\ast\nabla^{k-i}f\\
\end{eqnarray*}
One can see the expression (\ref{BIHk}) is homogeneous with respect
to covariant derivation. With out the symmetric condition, i.e.
$\nabla Rm=0$, there will be more terms occur, which is higher order
in $f$ but same order with respect to derivation. \qed

\lemma \label{HEHEF}Suppose $G:(M,g)\longrightarrow(N,h)$ satisfies
the Harmonic - Einstein metric (\ref{HaMEl}), then we have the
coupled system for the full curvature tensor and harmonic map:
\begin{eqnarray}
\Delta\nabla^kRm&=&Rm\ast
\nabla^{k}Rm+\sum_{i=0}^{k+2}\nabla^{i+1}G\ast
\nabla^{k+3-i}G+\sum_{i=1}^{k-1}\nabla^iRm\ast \nabla^{k-i}Rm\label{ECTH}\\
\Delta\nabla^k
G&=&\sum_{i=0}^{k-1}\nabla^iRm\ast\nabla^{k-i}G+\sum_{p=3}^{k+2}\sum_{i_1+\cdots+i_p=k-p+2}\nabla^{i_1+1}G\ast\cdots\ast\nabla^{i_p+1}G
\label{EHMH}
\end{eqnarray}

\proof First, let us recall that the full Riemannian curvature of
any Riemannian metric satisfies the following equation
\[\Delta Rm=\nabla^2Ric+Rm\ast Rm\]
and the well known Bochner formula on $M$, see \cite{Besse87}:
\begin{eqnarray*}
  \Delta\nabla^kRm &=&\nabla^k\Delta Rm+\sum_{i=0}^k\nabla^iRm\ast \nabla^{k-i}Rm \\
    &=& \nabla^{k+2}Ric+\sum_{i=0}^k\nabla^iRm\ast \nabla^{k-i}Rm
\end{eqnarray*}
Consequently, coupled with the Harmonic - Einstein equation
(\ref{HaMEl}), then Ricci curvature  is related to $dG$,
\[Ric=\lambda g+\langle dG,dG\rangle\]
Replacing $Ric$ term in the above expression, we have (\ref{ECTH}),
\begin{eqnarray*}
  \Delta\nabla^kRm &=& \nabla^{k+2}(dG\ast dG)+\sum_{i=0}^k\nabla^iRm\ast \nabla^{k-i}Rm\\
    &=& \sum_{i=0}^{k+2}\nabla^{i+1}G\ast \nabla^{k+3-i}G+\sum_{i=0}^k\nabla^iRm\ast \nabla^{k-i}Rm
\end{eqnarray*}
Since $G$ is a harmonic map, (\ref{EHMH}) follows from lemma
\ref{BTIF}. Moreover, we can drop the covariant derivatives of
$Rm^N$, since the target metric does not deform any more. \qed
\section{Local Regularity}
Now let us establish the $\epsilon$ - Regularity for Harmonic -
Einstein equation (\ref{HaMEl}) by divided the proof into several
lemmas.

\lemma \label{UBHF} Let $G:(M,g)\longrightarrow(N,h)$ satisfies the
Harmonic - Einstein equation (\ref{HaMEl}). Assume $(N,h)$ has
nonpositive sectional curvature and $(M,g)$ has bounded Sobolev
constant $C_S$. Then we have
\begin{equation}\label{BGGE}
\sup_{B_{\frac{R}{2}}(0)} |\nabla G|^2\leq
\frac{C}{R^n}\int_{B_R(0)}|\nabla G|^2.
\end{equation}
In other words, $|\nabla G|$ is bounded. Therefore, the $Ricci$
curvature is two sided bounded, and consequently, the volume of
geodesic ball is comparable with Euclidean Ball.

\proof Since $(N,h)$ has nonpositive sectional curvature, and
$\langle dG, dG\rangle$ is nonnegative, then by the Bochner formula
in lemma \ref{BTIF}, we obtain:
\begin{eqnarray*}
  \frac{1}{2}\Delta|dG|^2 &=& |\nabla \nabla f|^2-\sum_{\alpha,\beta}Rm^N(G_*e_\alpha,G_*e_\beta,G_*e_\alpha,G_*e_\beta)+\sum_{i}Ric^M(G^*\theta_i,G^*\theta_i)\\
 &\geq& \sum_{i}Ric^M(G^*\theta_i,G^*\theta_i)\\
 &\geq& \lambda |dG|^2
\end{eqnarray*}
With above equation in hand,  by elliptic Moser iteration  with
uniform Sobolev constant, the $L^\infty$ norm of $|dG|^2$ can be
bounded by $L^2$ of $|dG|^2$. Actually, this technique will be used
through our paper for more general tensors, so will not give the
detail here.  With more effort, one can also iterate by virtue of a
well known way, then $L^\infty$ norm can be bounded by $L^p$ norm of
$|dG|^2$ for any $p>0$. Now, we have $L^1$ - norm of $|dG|^2$, so we
have,
\[\sup_{B_{\frac{1}{2}}(0)} |dG|^2\leq C\int_{B_1(0)}|dG|^2.\]
For more details, see \cite{GT83}, \cite{SY94}, \cite{Simon96}. \qed

From now on, we will denote $\gamma=\frac{n}{n-2}$ through out this
paper.  $\phi$ will be a cut off function with $supp\,\phi\subset
B(0,r)$, and $\phi\equiv 1$ on $B(0,\tau)$ with $|\nabla \phi|\leq
\frac{2}{r-\tau}$. The estimation below will be affected by
different choices of $r$ and $\tau$, therefore we will choose proper
cut off function with respect to our purpose.

\lemma\label{IterI}(Iteration I:
$\parallel\!T\!\parallel_{L^{p\gamma}}$ estimation from $\Delta T$
with $\parallel\!T\!\parallel_{L^{p}}$) Let $T$ be a tensor, then
\begin{equation}\label{SbcLE}
\{\int (\phi|T|^{\frac{p}{2}})^{2\gamma}\}^\frac{1}{\gamma}
  \leq C(\int|\nabla
  \phi|^2|T|^p+p\phi^2|T|^{p-2}\langle
T,-\Delta T\rangle)
\end{equation}
\proof First, let us do some basic calculation. With Kato inequality
$|\nabla|T||\leq|\nabla T|$, then
\begin{eqnarray*}
-|T|\Delta |T|&=&-\frac{1}{2}\Delta |T|^2+|\nabla|T||^2\\
&=&-\langle
T,\Delta T\rangle-(|\nabla T |^2-|\nabla|T||^2)\\
&\leq& -\langle T,\Delta T\rangle
\end{eqnarray*}
Consequently, when $p\geq2$,
\begin{eqnarray*}
-\Delta |T|^{\frac{p}{2}} &=&-\frac{p}{2}|T|^{\frac{p}{2}-1}\Delta |T|-\frac{p}{2}(\frac{p}{2}-1)|T|^{\frac{p}{2}-2}|\nabla|T||^2 \\
&\leq& -\frac{p}{2}|T|^{\frac{p}{2}-2}\langle T,\Delta T\rangle
\end{eqnarray*}
Moreover, from the Schwarz inequality we obtain
\begin{eqnarray*}
\phi^2|\nabla|T|^{\frac{p}{2}}|^2&=&div(\phi^2|T|^{\frac{p}{2}}\nabla|T|^{\frac{p}{2}})-\phi^2|T|^{\frac{p}{2}}\Delta|T|^{\frac{p}{2}}-2\langle\phi\nabla|T|^{\frac{p}{2}},\nabla\phi|T|^{\frac{p}{2}}\rangle\\
&\leq&div(\phi^2|T|^{\frac{p}{2}}\nabla|T|^{\frac{p}{2}})-\frac{p}{2}\phi^2|T|^{p-2}\langle
T,\Delta
T\rangle+\delta\phi^2|\nabla|T|^{\frac{p}{2}}|^2+\frac{1}{\delta}|\nabla\phi|^2|T|^{p}
\end{eqnarray*}
Taking $\delta= \frac{1}{2}$, then the term
$\delta\phi^2|\nabla|T|^{\frac{p}{2}}|^2$ can be absorbed by the
left,
\begin{eqnarray*}
\phi^2|\nabla|T|^{\frac{p}{2}}|^2&\leq&2div(\phi^2|T|^{\frac{p}{2}}\nabla|T|^{\frac{p}{2}})-p\phi^2|T|^{p-2}\langle
T,\Delta T\rangle+4|\nabla\phi|^2|T|^{p}\\
\end{eqnarray*}
By Sobolev inequality,
\begin{eqnarray*}
  \{\int (\phi|T|^{\frac{p}{2}})^{2\gamma}\}^\frac{1}{\gamma} &\leq& C_S\int |\nabla(\phi |T|^{\frac{p}{2}})|^2 \\
  &\leq& C(\int|\nabla
  \phi|^2|T|^p+\int\phi^2|\nabla|T|^{\frac{p}{2}}|^2)\\
  &\leq& C(\int|\nabla
  \phi|^2|T|^p+p\phi^2|T|^{p-2}\langle
T,-\Delta T\rangle)
\end{eqnarray*}
Actually, this is nothing but the Moser iteration relation which is
generalized to tensor. \qed

\lemma \label{L2E} ($\parallel\!\nabla T\!\parallel_{L^2}$
estimation from $\Delta T$) Let $T$ be a tensor, then
\begin{equation}\label{GL2}
\int \phi^2|\nabla T|^2
  \leq C(\int\phi^2\langle T, -\Delta
T\rangle+\int|\nabla\phi|^2|T|^2)
\end{equation}
In particular, if
\[\Delta T=Rm\ast T+cT+\nabla X+Y\]
where $c$ is some constant, $X$, $Y$ are tensors. Then there exists
a constant $\epsilon=\epsilon(C_S,n)$, if
\[\{\int_{B(0,1)}
|Rm|^{\frac{n}{2}}\}^{\frac{n}{2}}\leq \epsilon,\] then we have
\begin{equation}\label{GWEL}
\int \phi^2|\nabla T|^2
  \leq C(\int\phi^2|X|^2+\int\phi^2|Y|^2+\int(\phi^2+|\nabla\phi|^2)|T|^2)
\end{equation}
\proof As did in lemma \ref{IterI}, we have,
\begin{eqnarray*}
\phi^2|\nabla T|^2&=&div(\phi^2\langle T, \nabla T\rangle-\phi^2\langle T, \Delta T\rangle-2\langle\phi\nabla T,\nabla\phi T\rangle\\
&\leq&div(\phi^2\langle T, \nabla T\rangle+\phi^2\langle T,-\Delta
T\rangle+\delta\phi^2|\nabla
T|^2+\frac{1}{\delta}|\nabla\phi|^2|T|^2
\end{eqnarray*}
Therefore
\[\int \phi^2|\nabla T|^2
  \leq C(\int\phi^2\langle T, -\Delta
T\rangle+\int|\nabla\phi|^2|T|^2)\] In particular, if
\[\Delta T=Rm\ast T+cT+\nabla X+Y\]
then the Laplacian term can be reduced to:
\begin{eqnarray*}
  \phi^2\langle T,-\Delta
T\rangle &=&\phi^2\langle T,-Rm\ast T-cT-\nabla X-Y\rangle \\
&=&-div(\phi^2\langle T,X\rangle)+\phi^2\langle\nabla
T,X\rangle+2\phi\langle\nabla\phi\otimes
T,X\rangle\\&&+\phi^2\langle T,-Rm\ast
T\rangle-\phi^2\langle T,Y\rangle-c\phi^2|T|^2\\
&\leq&-div(\phi^2\langle T,X\rangle)+\delta\phi^2|\nabla
T|^2+C(\phi^2|X|^2+|\nabla\phi|^2|T|^2\\&&+\phi^2|Rm||T|^2+\phi^2|T|^2+\phi^2|Y|^2)
\end{eqnarray*}
and hence
\begin{eqnarray*}
\phi^2|\nabla T|^2&\leq&div(\phi^2\langle T, \nabla T\rangle-div(\phi^2\langle T,X\rangle)\\
 &&+C(\phi^2|Rm||T|^2+\phi^2|X|^2+\phi^2|Y|^2+(\phi^2+|\nabla\phi|^2)|T|^2)
\end{eqnarray*}
By Sobolev inequality
\begin{eqnarray*}
  \{\int(\phi|T|)^{2\gamma}\}^\frac{1}{\gamma} &\leq&C(\int|\nabla \phi|^2|T|^2+\int\phi^2|\nabla T|^2)\\
&\leq&C(\int\phi^2|Rm||T|^2+\phi^2|X|^2+\phi^2|Y|^2+(\phi^2+|\nabla\phi|^2)|T|^2
\end{eqnarray*}
Moreover, if \[C\{\int_{B(0,1)}
|Rm|^{\frac{n}{2}}\}^{\frac{n}{2}}\leq \frac{1}{2},\] then the first
term
\[C\int\phi^2|Rm||T|^2\leq C\{\int|Rm|^{\frac{n}{2}}\}^{\frac{2}{n}}
\{\int(\phi|T|)^{2\gamma}\}^\frac{1}{\gamma}\leq \frac{1}{2}
\{\int(\phi|T|)^{2\gamma}\}^\frac{1}{\gamma}\] can be absorbed by
the left, thus we obtain
\begin{eqnarray*}
  \{\int(\phi|T|)^{2\gamma}\}^\frac{1}{\gamma}
  &\leq&C(\int\phi^2|X|^2+\int\phi^2|Y|^2+\int(\phi^2+|\nabla\phi|^2)|T|^2)
\end{eqnarray*}
which in turn implies
\begin{eqnarray*}
 \int \phi^2|\nabla T|^2 &\leq&C(\int\phi^2|X|^2+\phi^2|Y|^2+\int(\phi^2+|\nabla
\phi|^2)|T|^2)
\end{eqnarray*}
\qed

\lemma\label{ITLII}(Iteration II: $\parallel\!\nabla
T\!\parallel_{L^{p\gamma}}$ estimation from $\Delta\nabla T$ with
$\parallel\!\nabla T\!\parallel_{L^{p}}$)For any $p\geq2$,  there
exist a constant $\epsilon=\epsilon(C_S,n)$, if
\[\{\int_{B(0,1)}
|Rm|^{\frac{n}{2}}\}^{\frac{n}{2}}\leq \frac{\epsilon}{p},\]then
\begin{eqnarray}
\{\int (\phi|\nabla T|^{\frac{p}{2}})^{2\gamma}\}^\frac{1}{\gamma}
&\leq& C(\int\phi^2|\nabla T|^{p-2}|\Delta T|^2+\int\phi^2|\nabla
T|^{p-2}|Rm|^2|T|^2\nonumber\\&&+\int|\nabla\phi|^2|\nabla
T|^{p})\label{LpgLp}
\end{eqnarray}
In our case, $Rm=Rm^M+Rm^N(dG,dG)$, since $|\nabla G|^2|Rm^N|$ is
bounded, as explained in the proof, the above inequality holds if we
refer $Rm$ as $Rm^M$.

\proof Similarly, replacing $T$ by $\nabla T$ in lemma \ref{IterI},
we have
\begin{eqnarray*}
\phi^2|\nabla|\nabla T|^{\frac{p}{2}}|^2
&\leq&div(\phi^2|\nabla T|^{\frac{p}{2}}\nabla|\nabla
T|^{\frac{p}{2}})-p\phi^2|\nabla T|^{p-2}\langle \nabla
T,\Delta\nabla T\rangle+4|\nabla\phi|^2|\nabla T|^{p}
\end{eqnarray*}
By Bochner formula
\[\Delta\nabla T=\nabla\Delta T+\nabla (Rm\ast T)+Rm\ast \nabla T\]
and the we have
\begin{eqnarray*}
-p\phi^2|\nabla T|^{p-2}\langle \nabla T,\Delta\nabla T\rangle
&=&-p\phi^2|\nabla T|^{p-2}\langle \nabla
T,\nabla\Delta T+\nabla (Rm\ast T)+Rm\ast \nabla T\rangle\\
 &=&p \{-div(\phi^2|\nabla T|^{p-2}\langle
\nabla T,\Delta T+Rm\ast T\rangle)\\&&+\phi^2|\nabla T|^{p-2}\langle
\Delta T,\Delta T+Rm\ast T\rangle\\&& +(p-2)\phi^2|\nabla
T|^{p-3}\langle \nabla|\nabla T|\otimes\nabla T,\Delta T+Rm\ast
T\rangle\\&&+2\phi |\nabla T|^{p-2}\langle \nabla\phi\otimes\nabla
T,\Delta T+Rm\ast T\rangle\\&&+\phi^2|\nabla T|^{p-2}\langle \nabla
T,Rm\ast \nabla T\rangle\}
\end{eqnarray*}
For the term which involves second covariant derivative of $T$,
apply the Schwarz inequality,
\begin{eqnarray*}
&&p(p-2)\phi^2|\nabla T|^{p-3}\langle \nabla|\nabla T|\otimes\nabla
T,\Delta T+Rm\ast T\rangle\\
 &\leq& 2(p-2)\phi^2|\nabla
T|^{\frac{p-2}{2}}|\nabla|\nabla
T|^{\frac{p}{2}}|(|\Delta T|+|Rm||T|) \\
&\leq&\delta\phi^2|\nabla|\nabla
T|^{\frac{p}{2}}|^2+4p^2\phi^2|\nabla T|^{p-2}(|\Delta
T|^2+|Rm|^2|T|^2)
\end{eqnarray*}
Thus we obtain
\begin{eqnarray*}
\phi^2|\nabla|\nabla T|^{\frac{p}{2}}|^2&\leq&2div(\phi^2|\nabla
T|^{\frac{p}{2}}\nabla|\nabla T|^{\frac{p}{2}})-p\,div(\phi^2|\nabla
T|^{p-2}\langle \nabla T,\Delta T+Rm\ast
T\rangle)\\
&&+C(p^2\phi^2|\nabla T|^{p-2}|\Delta T|^2+p^2\phi^2|\nabla
T|^{p-2}|Rm|^2|T|^2\\&&+p\phi^2|Rm||\nabla
T|^{p}+p|\nabla\phi|^2|\nabla T|^{p})
\end{eqnarray*}
In our case, $Rm=Rm^M+Rm^N(dG,dG)$, since $|\nabla G|^2|Rm^N|$ is
bounded, so we have  $p\phi^2|(\nabla G)^2\ast Rm^N||\nabla
T|^{p}\leq Cp\phi^2|\nabla T|^{p}$, thus the above inequality holds
even if we refer $Rm$ as $Rm^M$.

Come back to the Sobolev inequality
\begin{eqnarray*}
  \{\int (\phi|\nabla T|^{\frac{p}{2}})^{2\gamma}\}^\frac{1}{\gamma}
  &\leq& C(\int|\nabla \phi|^2|\nabla T|^p+\int\phi^2|\nabla|\nabla
  T|^{\frac{p}{2}}|^2)\\
&\leq&C(p^2\int \phi^2|\nabla T|^{p-2}|\Delta
T|^2+p^2\int\phi^2|\nabla
T|^{p-2}|Rm|^2|T|^2\\&&+p\int\phi^2|Rm||\nabla
T|^{p}+p\int|\nabla\phi|^2|\nabla T|^{p})
\end{eqnarray*}
By H\"{o}lder inequality and our assumption of small integral of
curvature, the term
\[p\int\phi^2|Rm||\nabla T|^{p}\leq p\{\int|Rm|^{\frac{n}{2}}\}^{\frac{2}{n}} \{\int (\phi|\nabla T|^{\frac{p}{2}})^{2\gamma}\}^\frac{1}{\gamma}\]
can be absorbed by the left, thus we get
\begin{eqnarray*}
  \{\int (\phi|\nabla T|^{\frac{p}{2}})^{2\gamma}\}^\frac{1}{\gamma} &\leq&C(\int\phi^2|\nabla T|^{p-2}|\Delta T|^2+\int\phi^2|\nabla T|^{p-2}|Rm|^2|T|^2\\&&+\int|\nabla\phi|^2|\nabla
T|^{p})
\end{eqnarray*}
\qed

Now, we will use lemma \ref{IterI}, \ref{L2E} and \ref{ITLII} to get
some a priori estimation.

\lemma\label{G2L2} There exist a constant
$\epsilon=\epsilon(C_S,n)$, if
\[\{\int_{B(0,1)}
|Rm|^{\frac{n}{2}}\}^{\frac{n}{2}}\leq \epsilon,\] then
\begin{equation}
   \{\int \phi^2|\nabla^2G|^{2}\}^{\frac{1}{2}} \leq C\{\int|\nabla G|^2\}^{\frac{1}{2}}
\end{equation}
\proof Apply lemma \ref{L2E} to the equation
\[\Delta\nabla G=Rm^N\ast(\nabla G)^3+Ric^M\ast\nabla
G\] In (\ref{GWEL}), $T=\nabla G, X=0,Y=Rm^N\ast(\nabla
G)^3\approx|\nabla G|$, so we obtain the lemma.\qed

\theorem \label{P0E} For any $p\geq2$, there exist a constant
$\epsilon=\epsilon(C_S,n)$ such that if \[\{\int_{B(0,1)}
|Rm|^{\frac{n}{2}}\}^{\frac{2}{n}}\leq\frac{\epsilon}{p},\] then
\begin{eqnarray*}
  \{\int_{B(0,\frac{1}{2})}|\nabla^{2}G|^p\}^{\frac{1}{p}} &\leq& C(\{\int_{B(0,1)}|Rm|^{\frac{n}{2}}\}^{\frac{2}{n}}+\{\int_{B(0,1)}|\nabla
  G|^2\}^{\frac{1}{2}})\\
  \{\int_{B(0,\frac{1}{2})}|Rm|^p\}^{\frac{1}{p}} &\leq& C(\{\int_{B(0,1)}|Rm|^{\frac{n}{2}}\}^{\frac{2}{n}}+\{\int_{B(0,1)}|\nabla
  G|^2\}^{\frac{1}{2}})
\end{eqnarray*}
where $C=C(C_S,p,n)$.

\proof Since the curvature equation is coupled with the harmonic map
equation in lemma \ref{HEHEF}, we will see that we have to control
the two term in the lemma simultaneously. Recall the Iteration lemma
II \ref{ITLII}, let $T=\nabla G$, then we have
\begin{equation*}
\{\int (\phi|\nabla^2G|^{\frac{p}{2}})^{2\gamma}\}^\frac{1}{\gamma}
\leq C(\int\phi^2|\nabla^2G|^{p-2}|\Delta\nabla
G|^2+\int\phi^2|\nabla^2G|^{p-2}|Rm|^2|\nabla
G|^2+\int|\nabla\phi|^2|\nabla^2G|^{p})
\end{equation*}
By the equation
\[\Delta\nabla G=Rm^N\ast(\nabla G)^3+Ric^M\ast\nabla
G\] then
\begin{eqnarray*} |\Delta\nabla G|^2&\leq&C(|Rm|^2+|\nabla
G|^2)
\end{eqnarray*}
Therefore
\begin{equation*}
\{\int (\phi|\nabla^2G|^{\frac{p}{2}})^{2\gamma}\}^\frac{1}{\gamma}
\leq
C(\int\phi^2|\nabla^2G|^{p-2}|Rm|^2+\int\phi^2|\nabla^2G|^{p-2}|\nabla
G|^2+\int|\nabla\phi|^2|\nabla^2G|^{p})
\end{equation*}
Now apply the h\"{o}lder inequality with dual index
$\frac{p-2}{p}+\frac{2}{p}=1$, and we obtain,
\begin{eqnarray*}
  \{\int (\phi|\nabla^2G|^{\frac{p}{2}})^{2\gamma}\}^\frac{1}{\gamma}
&\leq&
C(\int\phi^2|Rm|^p+\int(\phi^2+|\nabla\phi|^2)|\nabla^2G|^{p}+\int\phi^2|\nabla
G|^{p})
\end{eqnarray*}
Furthermore, taking the $p$ - th root of both side with a trivial
inequality \[(\sum_{i=1}^Na_i)^{\frac{1}{p}}\leq
N^{\frac{1}{p}}\sum_{i=1}^Na_i^{\frac{1}{p}}\] we get
\begin{eqnarray}
\{\int \phi^{2\gamma}|\nabla^2G|^{p\gamma}\}^{\frac{1}{p\gamma}}
 &\leq&C^{\frac{1}{p}}(\sup|\nabla\phi|^{\frac{2}{p}}\{\int_{supp \phi}|\nabla^2G|^{p}\}^{\frac{1}{p}}+\{\int_{supp \phi}|Rm|^p\}^{\frac{1}{p}}\nonumber\\&&+\{\int_{supp \phi}|\nabla G|^{p}\}^{\frac{1}{p}})\label{aIGW2p}
\end{eqnarray}
Therefore, we obtain a priori estimation of $\nabla^2G$ but involves
$Rm$.

Now we turn to the estimation on $Rm$. Recall computation in lemma
\ref{IterI}, replacing $T$ by $Rm$,
\begin{eqnarray*}
\phi^2|\nabla|Rm|^{\frac{p}{2}}|^2
&\leq&2div(\phi^2|Rm|^{\frac{p}{2}}\nabla|Rm|^{\frac{p}{2}})-p\phi^2|Rm|^{p-2}\langle
Rm,\Delta Rm\rangle+4|\nabla\phi|^2|Rm|^{p}
\end{eqnarray*}
Combine with the equation of the curvature
\[\Delta Rm=\nabla^2Ric+Rm\ast Rm\]
we have
\begin{eqnarray*}
-\phi^2|Rm|^{p-2}\langle Rm,\Delta Rm\rangle
&=&-\phi^2|Rm|^{p-2}\langle
Rm,\nabla^2Ric+Rm\ast Rm\rangle\\
&=&-div(\phi^2|Rm|^{p-2}\langle Rm,\nabla
Ric\rangle)+\phi^2|Rm|^{p-2}\langle\delta Rm,\nabla
Ric\rangle\\&&+\frac{2(p-2)}{p}\phi^2|Rm|^{\frac{p}{2}-2}\langle
Rm,\nabla|Rm|^{\frac{p}{2}}\otimes\nabla Ric\rangle)
\\&&+2\phi|Rm|^{p-2}\langle Rm,\nabla\phi\otimes\nabla
Ric\rangle-\phi^2|Rm|^{p-2}\langle Rm,Rm\ast Rm\rangle\\
&\leq&-div(\phi^2|Rm|^{p-2}\langle Rm,\nabla
Ric\rangle)+\frac{\delta}{p}\phi^2|\nabla|Rm|^{\frac{p}{2}}|^2\\&&+C(p\phi^2|Rm|^{p-2}|\nabla
Ric|^2+|\nabla\phi|^2|Rm|^{p}+\phi^2|Rm|^{p+1})
\end{eqnarray*}
and consequently,
\begin{eqnarray*}
\phi^2|\nabla|Rm|^{\frac{p}{2}}|^2
&\leq&2div(\phi^2|Rm|^{\frac{p}{2}}\nabla|Rm|^{\frac{p}{2}})-p\,div(\phi^2|Rm|^{p-2}\langle
Rm,\nabla Ric\rangle)\\&&+C(p^2\phi^2|Rm|^{p-2}|\nabla
Ric|^2+p\phi^2|Rm|^{p+1}+p|\nabla\phi|^2|Rm|^{p})\\
\end{eqnarray*}
by  the Sobolev Inequality, absorbing the term $\int
p\phi^2|Rm|^{p+1}$ by the left, then \[ \{\int
(\phi|Rm|^{\frac{p}{2}})^{2\gamma}\}^\frac{1}{\gamma}\leq
C(\int|\nabla \phi|^2||Rm|^p+\int\phi^2|Rm|^{p-2}|\nabla
  Ric|^2)\]
On the other hand, coupled with the Harmonic - Einstein equation
(\ref{HaMEl}), $Ric=\lambda g+\langle dG,dG\rangle$. Therefore,
$\nabla Ric=\nabla^2G\ast\nabla G$, and then
\begin{eqnarray*}
  \{\int (\phi|Rm|^{\frac{p}{2}})^{2\gamma}\}^\frac{1}{\gamma}&\leq&C(\int|\nabla \phi|^2||Rm|^p+\int\phi^2|Rm|^{p-2}|\nabla^2G|^2)\\
  &\leq&C(\int(\phi^2+|\nabla \phi|^2)||Rm|^p+\int
  \phi^2|\nabla^2G|^p)
\end{eqnarray*}
As did for $\nabla^2G$ (\ref{aIGW2p}), we have
\begin{eqnarray}
  \{\int\phi^{2\gamma}|Rm|^{p\gamma}\}^\frac{1}{p\gamma}&\leq&C^{\frac{1}{p}}(\sup|\nabla
  \phi|^{\frac{2}{p}}\{\int_{supp\;\phi}|Rm|^p\}^{\frac{1}{p}}+\{\int_{supp\;\phi}|\nabla^2G|^p\}^{\frac{1}{p}})\label{aIRmp}
\end{eqnarray}
By taking $p_i=2\gamma^i$, $supp\,\phi_i \subset B_i:=B(0,
\frac{1}{2}+(\frac{1}{2})^{i})$, $\phi_i\equiv 1$ on $B_{i+1}$, and
$|\nabla \phi_i|\leq 2^{i+2}$. Define
\begin{eqnarray*}
\Phi_i(\nabla^2G)&=&\{\int_{B_i}
|\nabla^2G|^{2\gamma^i}\}^{\frac{1}{2\gamma^i}} \\
\Psi_i(Rm)&=&\{\int_{B_i} |Rm|^{2\gamma^i}\}^{\frac{1}{2\gamma^i}}
\end{eqnarray*}
With (\ref{aIGW2p}) and (\ref{aIRmp}), we obtain a coupled iteration
relation
\begin{equation*}
  \begin{bmatrix}
    \Phi_{i+1}(\nabla^2G) \\
    \Psi_{i+1}(Rm) \\
  \end{bmatrix}
\leq C^{\frac{1}{2\gamma^i}}
  \begin{bmatrix}
    2^{\frac{i}{2\gamma^i}} & 1 \\
    1 &2^{\frac{i}{2\gamma^i}} \\
  \end{bmatrix}
 \begin{bmatrix}
   \Phi_i(\nabla^2G) \\
 \Psi_i(Rm) \\
  \end{bmatrix}
 + C^{\frac{1}{2\gamma^i}}
 \begin{bmatrix}
    \Phi_i(\nabla G) \\
0 \\
  \end{bmatrix}
\end{equation*}
Denote $\lambda_i=2^{\frac{i}{2\gamma^i}}$, and use the following
fact on matrix:
\begin{equation*}
\begin{bmatrix}
    \lambda & 1 \\
    1 &  \lambda \\
  \end{bmatrix}= \frac{1}{\sqrt{2}}
\begin{bmatrix}
    1 & -1 \\
    1 & 1 \\
   \end{bmatrix}
\begin{bmatrix}
    \lambda+1&0\\
   0 &  \lambda-1 \\
 \end{bmatrix}\frac{1}{\sqrt{2}}
\begin{bmatrix}
    1 & -1 \\
    1 & 1 \\
 \end{bmatrix}^T
\end{equation*}
By iterating on $i$, and consequently,
\begin{eqnarray}
&&\begin{bmatrix}
    \Phi_{i+1}(\nabla^2G) \\
    \Psi_{i+1}(Rm) \\
\end{bmatrix}\nonumber\\&\leq& \frac{1}{2}C^{\sum_{j=0}^i\frac{1}{2\gamma^j}}
\begin{bmatrix}
    \Pi_{j=0}^i(\lambda_j+1)+\Pi_{j=0}^i(\lambda_j-1)&\Pi_{j=0}^i(\lambda_j+1)-\Pi_{j=0}^i(\lambda_j-1)\\
   \Pi_{j=0}^i(\lambda_j+1)-\Pi_{j=0}^i(\lambda_j-1) & \Pi_{j=0}^i(\lambda_i+1)+\Pi_{j=0}^i(\lambda_j-1) \\
\end{bmatrix}
\begin{bmatrix}
\Phi_0(\nabla^2G) \\
\Psi_0(Rm) \\
\end{bmatrix}\nonumber\\
&&+
\sum_{j=0}^iC^{\sum_{k=j}^i\frac{1}{2\gamma^k}}2^{\sum_{k=j+1}^{i}\frac{k}{2\gamma^k}}
\begin{bmatrix}
\Phi_j(\nabla G) \\
0 \\
\end{bmatrix}\nonumber\\
&\leq& C\begin{bmatrix}
            e^{c_1(n)i}+c_2(n) & e^{c_1(n)i}-c_2(n) \\
            e^{c_1(n)i}-c_2(n) & e^{c_1(n)i}+c_2(n) \\
\end{bmatrix}
\begin{bmatrix}
\Phi_0(\nabla^2G) \\
\Psi_0(Rm) \\
\end{bmatrix}+ Ci
\begin{bmatrix}
\Phi_0(\nabla G) \\
0 \\
\end{bmatrix}\label{ICST}
\end{eqnarray}
where we have used the facts:
$\Pi_{j=0}^i(2^{\frac{j}{2\gamma^j}}+1)=e^{\sum_{j=0}^i\ln(2^{\frac{j}{2\gamma^j}}+1)}\leq
e^{\sum_{j=0}^i2^{\frac{j}{2\gamma^j}}}\leq e^{c_1(n)i}$,
$\Pi_{j=0}^i(2^{\frac{j}{2\gamma^j}}-1)\leq c_2(n)$, $\sum
\frac{1}{q\gamma^j}=\frac{\gamma}{q(\gamma-1)}$, $\sum
\frac{i}{q\gamma^j}=\frac{\gamma}{q(\gamma-1)^2}$.

For the initial condition, lemma \ref{G2L2} implies
\[\Phi_{0}(\nabla^2G)=\{\int_{B_0}
|\nabla^2G|^{2}\}^{\frac{1}{2}}\leq C\{\int_{B(0,1)}|\nabla
G|^{2}\}^{\frac{1}{2}}\] and also as explained in lemma \ref{UBHF},
the volume of geodesic ball is comparable with Euclidean ball, then
\[\Psi_{0}(Rm)=\{\int_{B_0} |Rm|^{2}\}^{\frac{1}{2}}\leq
C\{\int_{B(0,1)} |Rm|^{\frac{n}{2}}\}^{\frac{2}{n}}\]

Finally, come back to (\ref{ICST}), we  obtain
\begin{eqnarray*}
\{\int_{B_{i}}|\nabla^2G|^{2\gamma^{i}}\}^{\frac{1}{2\gamma^{i}}}&=&C(\{\int|Rm|^{\frac{n}{2}}\}^{\frac{2}{n}}+\{\int|\nabla
G|^{2}\}^{\frac{1}{2}})\\
\{\int_{B_{i}}|Rm|^{2\gamma^{i}}\}^\frac{1}{2\gamma^{i}}&=&C(\{\int|Rm|^{\frac{n}{2}}\}^{\frac{2}{n}}+\{\int|\nabla
G|^{2}\}^{\frac{1}{2}})
\end{eqnarray*}
where $C=C(n,C_S, i)$. \qed

Before we starting to prove the $\epsilon$ - regularity \ref{eRegT},
let us state the Moser iteration lemma. \lemma (Moser iteration
\cite{GT83} \cite{BKN89})\label{MIL} Suppose a nonnegative function
$u$ satisfies $\Delta u\geq - f u- h - cu$ with $f\in L^{q}, q
> \frac{n}{2}$, $g\in L^{q'}, q'> \frac{n}{2}$, $c$ is some constant, and $u\in L^p$ for some $p\in [p_0,
p_1]$ where $p_0 > 1$. Since we do analysis  on manifolds, we also
assume bounded  $C_S$  and Euclidean volume growth, i.e.
$vol(B(0,r))\leq Vr^n$. Then there exists a constant $C = C(p_0,p_1,
C_S, V, c, \parallel\!f\!\parallel_{L^q})$ so that
\begin{equation}
 \sup_{B(0,\frac{r}{2})}|u| \leq C r^{-\frac{n}{p}}\{\int_{B(0,r)}
 |u|^p\}^\frac{1}{p}+Cr^{-\frac{n}{q'}}\{\int_{B(0,r)}|h|^{q'}\}^{\frac{1}{q'}}
\end{equation}

Since all the inequalities in the main theorem \ref{eRegT} are scale
invariant, we may assume $r=1$ for simplicity, and then theorem
\ref{epsLe} is equivalent to theorem \ref{eRegT}.
 \theorem\label{GRmWkp}For any
$k\in\mathbb{N}$ and $p\geq 2$, there exist a constant
$\epsilon=\epsilon(C_S,n)$ such that if
\[\{\int_{B(0,1)}
|Rm|^{\frac{n}{2}}\}^{\frac{2}{n}}\leq\frac{\epsilon}{p},\] then
\begin{eqnarray}
  \{\int_{B(0,\frac{1}{2})}|\nabla^{k+2}G|^2\}^{\frac{1}{2}}
  &\leq&C(\{\int_{B(0,1)}|Rm|^{\frac{n}{2}}\}^{\frac{2}{n}}+\{\int_{B(0,1)}|\nabla
  G|^2\}^{\frac{1}{2}})\label{GWk22}\\
  \{\int_{B(0,\frac{1}{2})}|\nabla^{k+2}G|^p\}^{\frac{1}{p}} &\leq&C(\{\int_{B(0,1)}|Rm|^{\frac{n}{2}}\}^{\frac{2}{n}}+\{\int_{B(0,1)}|\nabla
  G|^2\}^{\frac{1}{2}})\label{GWk2p}\\
   \{\int_{B(0,\frac{1}{2})}|\nabla^kRm|^2\}^{\frac{1}{2}} &\leq&C(\{\int_{B(0,1)}|Rm|^{\frac{n}{2}}\}^{\frac{2}{n}}+\{\int_{B(0,1)}|\nabla
  G|^2\}^{\frac{1}{2}})\label{RmWk2}\\
  \{\int_{B(0,\frac{1}{2})}|\nabla^kRm|^p\}^{\frac{1}{p}} &\leq&C(\{\int_{B(0,1)}|Rm|^{\frac{n}{2}}\}^{\frac{2}{n}}+\{\int_{B(0,1)}|\nabla
  G|^2\}^{\frac{1}{2}})\label{RmWkp}
\end{eqnarray}
where $C=C(C_S,k,p,n)$.

\theorem\label{epsLe}($\epsilon$ - regularity) There exist a
$\epsilon=\epsilon(C_S,n)$ such that if
\[\{\int_{B(0,1)}
|Rm|^{\frac{n}{2}}\}^{\frac{2}{n}}\leq\epsilon\] then for any
$k\in\mathbb{N}$, we have
\begin{eqnarray}
\sup_{B(0,\frac{1}{2})}|\nabla^{k+1}G|&\leq&
C(\{\int_{B(0,1)}|Rm|^{\frac{n}{2}}\}^{\frac{2}{n}}+\{\int_{B(0,1)}|\nabla
  G|^2\}^{\frac{1}{2}}) \label{GUBk2}\\
\sup_{B(0,\frac{1}{2})}|\nabla^{k-1}Rm|&\leq&
C(\{\int_{B(0,1)}|Rm|^{\frac{n}{2}}\}^{\frac{2}{n}}+\{\int_{B(0,1)}|\nabla
  G|^2\}^{\frac{1}{2}})\label{RmUBk}
\end{eqnarray}
where $C=C(C_S,k,n)$.

\remark We will proof theorem \ref{GRmWkp} and theorem \ref{epsLe}
together by induction on k. Strictly speaking, for bigger $k$ one
need shrink the ball further after each step of the iteration. We
will take this for granted for short, but note that the constant $C$
will depends on $k$.

Before we starting the proof, let us present the main idea. As did
in the case $k=0$, see theorem \ref{P0E}, if we apply (\ref{GL2}) on
the equation $\Delta \nabla^{k+1}G$, we can get the $L^2$ estimation
on $\nabla^{k+2}G$; Similarly, if we apply (\ref{GL2}) on the
equation $\Delta \nabla^{k-1}Rm$, we can get the $L^2$ estimation on
$\nabla^{k}Rm$, thus we get (\ref{GWk22}) and (\ref{RmWk2}). If
$2\gamma>\frac{n}{2}$, namely $n\leq 5$, with Sobolev inequality,
(\ref{GWk22}) and (\ref{RmWk2}) is enough for Moser iteration to
bound the curvature and harmonic map. While for the higher dimension
case, we must apply the iteration lemma II \ref{ITLII} to the
equation $\Delta \nabla^{k+2}G$ and $\Delta\nabla^kRm$ to improve
the integrality order up to $p>\frac{n}{2}$. However, as we have
already seen in the case $k=0$, we can not get a priori estimation
for $|\nabla^k Rm|$ and $|\nabla^{k+2}G|$ separately like $n\leq 5$,
but get (\ref{GWk2p}) and (\ref{RmWkp}) simultaneously,  since our
equation is a coupled system. Once the integrality order is bigger
than $p>\frac{n}{2}$, one can apply the Moser iteration lemma
\ref{MIL} to get the $L^\infty$ estimate for $|\nabla^{k+1}G|$  and
$|\nabla^{k-1} Rm|$ , therefore we get (\ref{GUBk2}) and
(\ref{RmUBk}).

\proof We have already proved the case $k=0$ in theorem \ref{P0E}.
Moreover, we will see the case $k=1$ in theorem \ref{GRmWkp} does
not require theorem \ref{epsLe}. The theorem \ref{epsLe} will begin
from $k=1$ and the case $k=1$ will be proved in step III, which
require the case $k=1$ in theorem \ref{GRmWkp}. Thus the induction
process is well ordered. Now we assume all the inequalities in
theorem \ref{epsLe} and \ref{GRmWkp} hold for the case from $k=0$
through out to $k-1$.

Step I: Recall lemma \ref{HEHEF}, for the Harmonic - Einstein
equation (\ref{HaMEl}), we have the coupled system (\ref{ECTH}) and
(\ref{EHMH}) for the full curvature tensor $Rm$ and harmonic map
$G$. Now if we apply lemma \ref{L2E} on the equation (\ref{EHMH})
\begin{eqnarray*}
\Delta\nabla^{k+1}G &=& \sum_{i=0}^{k}\nabla^iRm\ast\nabla^{k-i+1}G+\sum_{p=3}^{k+3}\sum_{i_1+\cdots+i_p=k-p+3}\nabla^{i_1+1}G\ast\cdots\ast\nabla^{i_p+1}G \\
&=& \nabla(\nabla^{k-1}Rm\ast\nabla G)+Rm\ast\nabla^{k+1}G+\sum_{i=1}^{k-1}\nabla^iRm\ast\nabla^{k-i+1}G\\
&&+\nabla G\ast\nabla
G\ast\nabla^{k+1}G+\sum_{p=3}^{k+3}\sum_{i_1+\cdots+i_p=k-p+3,
i_*<k}\nabla^{i_1+1}G\ast\cdots\ast\nabla^{i_p+1}G
\end{eqnarray*}
with $T=\nabla^{k+1}G$, $c=|\nabla G|^2$, $X=\nabla^{k-1}
Rm\ast\nabla G$. When $k=1$, then 
\[\Delta\nabla^2 G=\nabla(Rm\ast\nabla G)+Rm\ast\nabla^2G+(\nabla
G)^2\ast\nabla^2G+(\nabla G)^4,\] Therefore, $Y=(\nabla G)^4$ in our
notation,  and $|Y|\leq C\{\int_{B(0,1)}|\nabla
  G|^2\}^{\frac{1}{2}}$ without the induction in theorem
\ref{epsLe}; When $k\geq2$, by induction, (\ref{GUBk2}) and
(\ref{RmUBk}) hold up to $k-1$, namely, $|\nabla^{j}Rm|$ and
$|\nabla^{j+2}G|$ are bounded for $j\leq k-2$, then
\begin{eqnarray*}
|Y|&=&|\sum_{i=1}^{k-1}\nabla^iRm\ast\nabla^{k-i+1}G+\sum_{p=3}^{k+3}\sum_{i_1+\cdots+i_p=k-p+3,
i_*<k}\nabla^{i_1+1}G\ast\cdots\ast\nabla^{i_p+1}G|\\
&\leq&
C|\nabla^{k-1}Rm|+C(\{\int_{B(0,1)}|Rm|^{\frac{n}{2}}\}^{\frac{2}{n}}+\{\int_{B(0,1)}|\nabla
  G|^2\}^{\frac{1}{2}})
\end{eqnarray*}
Come back to lemma \ref{L2E}, we have
\begin{eqnarray*}
\int\phi^2|\nabla^{k+2}G|^2&\leq&C(\int\phi^2|X|^2+\int(\phi^2+|\nabla\phi|^2)|T|^2+\int\phi^2|Y|^2) \\
&\leq&C(\int\phi^2|\nabla^{k-1}Rm|^2+\int(\phi^2+|\nabla\phi|^2)|\nabla^{k+1}G|^2\\
&&+(\{\int_{B(0,1)}|Rm|^{\frac{n}{2}}\}^{\frac{2}{n}}+\{\int_{B(0,1)}|\nabla
  G|^2\}^{\frac{1}{2}})^2)
\end{eqnarray*}
By induction, (\ref{GWk22}) and (\ref{RmWk2}) hold for with $k-1$,
so (\ref{GWk22}) holds for $k$.

Similarly, if we apply (\ref{GL2}) on the equation \ref{ECTH},
\begin{equation*}
    \Delta\nabla^{k-1}Rm=Rm\ast \nabla^{k-1}Rm+\sum_{i=0}^{k+1}\nabla^{i+1}G\ast \nabla^{k+2-i}G+\sum_{i=1}^{k-2}\nabla^iRm\ast \nabla^{k-i-1}Rm
\end{equation*}
with $T=\nabla^{k-1}Rm$, $c=0$, $X=0$, and
\begin{eqnarray*}
|Y|&=&|\nabla G\ast \nabla^{k+2}G+\sum_{i=1}^{k}\nabla^{i+1}G\ast
\nabla^{k+2-i}G+\sum_{i=1}^{k-2}\nabla^iRm\ast \nabla^{k-i-1}Rm| \\
&\leq &
(|\nabla^{k+2}G|+\sum_{i=1}^{k}|\nabla^{i+1}G|^2+\sum_{i=1}^{k-2}|\nabla^iRm|^2)
\end{eqnarray*}
then we have
\begin{eqnarray*}
\int\phi^2|\nabla^{k}Rm|^2&\leq&C(\int(\phi^2+|\nabla\phi|^2)|T|^2+\int\phi^2|Y|^2) \\
&=&C(\int(\phi^2+|\nabla\phi|^2)|\nabla^{k-1}
Rm|^2+\int\phi^2|\nabla^{k+2}G|^2\\&&
+\sum_{i=1}^{k}\int\phi^2|\nabla^{i+1}G|^4+\sum_{i=1}^{k-2}\int\phi^2|\nabla^iRm|^4)
\end{eqnarray*}
By induction, (\ref{GWk22}) - (\ref{RmWkp}) hold up to $k-1$, and
(\ref{GWk22}) holds for $k$ which is proved just now, then we have
(\ref{RmWk2}) for $k$.

Step II, apply the iteration lemma II \ref{ITLII} to  the equation
(\ref{EHMH}), with $T=\nabla^{k+1}G$, then we have
\begin{eqnarray*}
\{\int
(\phi|\nabla^{k+2}G|^{\frac{p}{2}})^{2\gamma}\}^\frac{1}{\gamma}
&\leq& C(\int\phi^2|\nabla^{k+2}G|^{p-2}|\Delta\nabla^{k+1}
G|^2\\
&&+\int\phi^2|\nabla^{k+2}G|^{p-2}|Rm|^2|\nabla^{k+1}
G|^2+\int|\nabla\phi|^2||\nabla^{k+2}G|^{p})
\end{eqnarray*}
Applying  Schwartz inequality to the equation (\ref{EHMH}),
$\Delta\nabla^{k+1}G$: when $k=1$, we have
\begin{eqnarray*}
|\Delta\nabla^2 G|^2 &\leq&C(|Rm|^2|\nabla^2 G|^2+|\nabla
Rm|^2+|\nabla^2 G|^2+\int_{B(0,1)}|\nabla
  G|^2)
\end{eqnarray*}
therefore, we do not use the induction in theorem \ref{epsLe}; when
$k\geq2$, by induction, (\ref{GUBk2}) and (\ref{RmUBk}) hold up to
$k-1$, then we have
\begin{eqnarray*}
|\Delta\nabla^{k+1}G|^2 &\leq&C(|Rm|^2|\nabla^{k+1}G|^2+|\nabla^{k}
Rm|^2+|\nabla^{k-1}
Rm|^2+|\nabla^{k+1}G|^2\\&&+(\{\int_{B(0,1)}|Rm|^{\frac{n}{2}}\}^{\frac{2}{n}}+\{\int_{B(0,1)}|\nabla
  G|^2\}^{\frac{1}{2}})^2)
\end{eqnarray*}
Replacing the Laplacian term in the above integral inequality,
\begin{eqnarray*}
\{\int
(\phi|\nabla^{k+2}G|^{\frac{p}{2}})^{2\gamma}\}^\frac{1}{\gamma}
&\leq& C(\int\phi^2|\nabla^{k+2}G|^{p-2}|Rm|^2|\nabla^{k+1}
G|^2+\int\phi^2|\nabla^{k+2}G|^{p-2}|\nabla^{k}
Rm|^2\\
&&+\int\phi^2|\nabla^{k+2}G|^{p-2}(|\nabla^{k-1}
Rm|^2+|\nabla^{k+1}G|^2)+\int|\nabla\phi|^2||\nabla^{k+2}G|^{p}\\&&+\int\phi^2||\nabla^{k+2}G|^{p-2}(\{\int_{B(0,1)}|Rm|^{\frac{n}{2}}\}^{\frac{2}{n}}+\{\int_{B(0,1)}|\nabla
  G|^2\}^{\frac{1}{2}})^2)\end{eqnarray*}
By H\"{o}lder inequality with $\frac{p-2}{p}+\frac{2}{p}=1$, we have
\begin{eqnarray}
\{\int
(\phi|\nabla^{k+2}G|^{\frac{p}{2}})^{2\gamma}\}^\frac{1}{\gamma}
&\leq&
C(\int(\phi^2+|\nabla\phi|^2)|\nabla^{k+2}G|^{p}+\int\phi^2|\nabla^{k}
Rm|^p+\int\phi^2|\nabla^{k-1}Rm|^{p}\nonumber\\
&&+\int\phi^2|\nabla^{k+1}G|^{p}+\int\phi^2|\nabla^{k+1}G|^{2p}+\int\phi^2|Rm|^{2p}\nonumber\\
&&+(\{\int_{B(0,1)}|Rm|^{\frac{n}{2}}\}^{\frac{2}{n}}+\{\int_{B(0,1)}|\nabla
  G|^2\}^{\frac{1}{2}})^p)\nonumber\\
  &\leq&
C(\int(\phi^2+|\nabla\phi|^2)|\nabla^{k+2}G|^{p}+\int\phi^2|\nabla^{k}
Rm|^p\nonumber\\
&&+(\{\int_{B(0,1)}|Rm|^{\frac{n}{2}}\}^{\frac{2}{n}}+\{\int_{B(0,1)}|\nabla
  G|^2\}^{\frac{1}{2}})^p)\label{aIGWk2p}
\end{eqnarray}
where the last inequality follows from the induction, namely,
(\ref{GWk22}) - (\ref{RmWkp}) for  $k-1$. Therefore we get a priori
estimation on $|\nabla^{k+2}G|$ but involves $\nabla^kRm$.

Now we turn to the estimation of $\nabla^kRm$. By the Iteration
lemma II \ref{ITLII} again to (\ref{ECTH}),  let $T=\nabla^{k-1}Rm$,
\begin{eqnarray*}
  \{\int (\phi|\nabla^kRm|^{\frac{p}{2}})^{2\gamma}\}^{\frac{1}{\gamma}} &\leq&C(\int \phi^{2}|\nabla^kRm|^{p-2}|\Delta\nabla^{k-1}Rm|^2\\&&+\int \phi^{2}|\nabla^kRm|^{p-2}|Rm|^2|\nabla^{k-1}Rm|^2+\int|\nabla \phi|^2|\nabla^kRm|^{p})
\end{eqnarray*}
Applying  Schwartz inequality to (\ref{ECTH}),
\begin{eqnarray*}
|\Delta\nabla^{k-1}Rm|^2&\leq&C(|\nabla^{k+2}G|^2+\sum_{i=1}^{k}|\nabla^{i+1}G|^2|\nabla^{k+2-i}G|^2+\sum_{i=0}^{k-1}|\nabla^iRm|^2|\nabla^{k-i-1}Rm|^2)
\end{eqnarray*}
Then we have
\begin{eqnarray*}
  &&\{\int (\phi|\nabla^kRm|^{\frac{p}{2}})^{2\gamma}\}^{\frac{1}{\gamma}}\\
   &\leq&C(\int \phi^{2}|\nabla^kRm|^{p-2}|\nabla^{k+2}G|^2+\sum_{i=1}^{k}\int \phi^{2}|\nabla^kRm|^{p-2}|\nabla^{i+1}G|^2|\nabla^{k+2-i}G|^2\\
  &&+\sum_{i=0}^{k-1}\int \phi^{2}|\nabla^kRm|^{p-2}|\nabla^iRm|^2|\nabla^{k-i-1}Rm|^2+\int|\nabla
  \phi|^2|\nabla^kRm|^{p})
\end{eqnarray*}
By H\"{o}lder inequality with $\frac{p-2}{p}+\frac{2}{p}=1$,
\begin{eqnarray}
  \{\int (\phi|\nabla^kRm|^{\frac{p}{2}})^{2\gamma}\}^{\frac{1}{\gamma}}
  &\leq&C(\int\phi^{2}|\nabla^{k+2}G|^p+\int(\phi^2+|\nabla
  \phi|^2)|\nabla^kRm|^{p}\nonumber\\
  &&+\sum_{i=1}^{k}\int \phi^{2}|\nabla^{i+1}G|^{2p}+\sum_{i=0}^{k-1}\int\phi^{2}|\nabla^iRm|^{2p})\nonumber\\
  &\leq&C(\int\phi^{2}|\nabla^{k+2}G|^p+\int(\phi^2+|\nabla
  \phi|^2)|\nabla^kRm|^{p}\nonumber\\
  &&+(\{\int_{B(0,1)}|Rm|^{\frac{n}{2}}\}^{\frac{2}{n}}+\{\int_{B(0,1)}|\nabla
  G|^2\}^{\frac{1}{2}})^p)\label{aIRmkp}
\end{eqnarray}
As did in theorem \ref{P0E}, see (\ref{ICST}), taking
$p_i=2\gamma^i$, $supp\,\phi_i \subset B_i:=B(0,
\frac{1}{2}+(\frac{1}{2})^{i})$ and $\phi_i\equiv 1$ on $B_{i+1}$,
and $|\nabla \phi|\leq 2^{i+2}$. Define
\begin{eqnarray*}
\Phi_i(\nabla^{k+2}G)&=&\{\int_{B_i}
|\nabla^{k+2}G|^{2\gamma^i}\}^{\frac{1}{2\gamma^i}} \\
\Psi_i(\nabla^kRm)&=&\{\int_{B_i}
|\nabla^kRm|^{2\gamma^i}\}^{\frac{1}{2\gamma^i}}
\end{eqnarray*}
and
\[C_0=C(\{\int|Rm|^{\frac{n}{2}}\}^{\frac{2}{n}}+\{\int|\nabla
  G|^2\}^{\frac{1}{2}})\]
With (\ref{aIGWk2p}) and (\ref{aIRmkp}), we obtain a coupled
iteration sequence
\begin{equation*}
  \begin{bmatrix}
    \Phi_{i+1}(\nabla^{k+2}G) \\
    \Psi_{i+1}(\nabla^kRm) \\
  \end{bmatrix}
\leq C^{\frac{1}{2\gamma^i}}
  \begin{bmatrix}
    2^{\frac{i}{2\gamma^i}} & 1 \\
    1 &2^{\frac{i}{2\gamma^i}} \\
  \end{bmatrix}
 \begin{bmatrix}
   \Phi_i(\nabla^{k+2}G) \\
 \Psi_i(\nabla^kRm) \\
  \end{bmatrix}
 + C^{\frac{1}{2\gamma^i}}
 \begin{bmatrix}
    C_0 \\
C_0\\
  \end{bmatrix}
\end{equation*}
 With the same iteration process as in
(\ref{ICST}), we obtain
\begin{equation*}
  \begin{bmatrix}
    \Phi_{i+1}(\nabla^{k+2}G) \\
    \Psi_{i+1}(\nabla^kRm) \\
  \end{bmatrix}
\leq C
\begin{bmatrix}
    e^{c_1(n)i}+c_2(n) & e^{c_1(n)i}-c_2(n) \\
    e^{c_1(n)i}-c_2(n) & e^{c_1(n)i}+c_2(n)\\
  \end{bmatrix}
   \begin{bmatrix}
 \Phi_{0}(\nabla^{k+2}G) \\
\Psi_{0}(\nabla^kRm)\\
  \end{bmatrix}
 + Ci
 \begin{bmatrix}
    C_0 \\
C_0\\
\end{bmatrix}
\end{equation*}
The initial condition (\ref{GWk22}) and (\ref{RmWk2}) for $i=0$, are
proved in step I, and therefore we proved (\ref{GWk2p}) and
(\ref{RmWkp}) for the case $k$.

Step III: We will apply Moser iteration to get the $L^\infty$
estimation (\ref{GUBk2}) and (\ref{RmUBk}).

For initial case $k=1$, once we have $L^{q} (q>\frac{n}{2})$ bound
of $\nabla^2 Ric$ and $L^{q'} (q'>\frac{n}{2})$ bound of $Rm$, by
lemma \ref{MIL}, we can apply the Moser iteration  to the equation
    \[-\Delta |Rm|\leq C(n)|\nabla^2 Ric|+ C(n)|Rm||Rm|\]
to obtain the $L^\infty$ estimation of the full curvature tensor.
Note that in the proof of step I and II for the case $k=1$, we do
not need the induction in theorem \ref{epsLe}, therefore we have
(\ref{GWk22}) - (\ref{RmWkp}) hold for $k=1$. If we take $p=2^i>n$,
namely, $i=\lfloor\frac{\ln n}{\ln 2}\rfloor+1$ in theorem
\ref{GRmWkp}, then we have $L^{\frac{p}{2}}$ bound of
$\sum_{j=1}^3|\nabla^jG|^2\geq|\nabla^2 Ric|$ and $L^{p}$ bound of
$|Rm|$, which implies the $L^\infty$ bound of the full curvature
tensor. On the other hand, we also have $L^p (p>n)$ bound of $\nabla
Rm$, the same argument on the equation
 \[-\Delta|\nabla^2G|\leq C(n)(|Rm||\nabla^2G|+|\nabla G|^2|\nabla^2G|+|\nabla G||\nabla Rm|^2+|\nabla G|^4)\]
will give the $L^\infty$ estimation on the derivation of $G$ up to
second order.

For any $k>1$, we have assumed,  by induction, (\ref{GUBk2}) and
(\ref{RmUBk}) hold up to $k-1$.

Apply Moser iteration lemma \ref{MIL} to the equation,
\begin{eqnarray*}
-\Delta|\nabla^{k+1}G|&\leq&C(n)(|Rm||\nabla^{k+1}G|+|\nabla
G|^2|\nabla^{k+1}G|+|\nabla
G||\nabla^{k}Rm|+|\nabla^{2}G||\nabla^{k-1}Rm|)\\&&+C(\{\int|Rm|^{\frac{n}{2}}\}^{\frac{2}{n}}+\{\int|\nabla
  G|^2\}^{\frac{1}{2}})
\end{eqnarray*}
we have $L^{p} (p>\frac{n}{2})$ norm of $|\nabla
G||\nabla^{k}Rm|+|\nabla^{2}G||\nabla^{k-1}Rm|+C(\{\int|Rm|^{\frac{n}{2}}\}^{\frac{2}{n}}+\{\int|\nabla
  G|^2\}^{\frac{1}{2}})$,
and also $L^p$ norm of $|\nabla^{k+1}G|$ by (\ref{GWk2p}) -
(\ref{RmWkp}). Apply lemma \ref{MIL} once more, we obtain the
$L^\infty$ estimation of $|\nabla^{k+1}G|$.

Similarly, on the equation
\begin{equation*}
   \Delta\nabla^{k-1}Rm=Rm\ast \nabla^{k-1}Rm+\sum_{i=0}^{k+1}\nabla^{i+1}G\ast \nabla^{k+2-i}G+\sum_{i=1}^{k-2}\nabla^iRm\ast \nabla^{k-i-1}Rm
\end{equation*}
by induction on (\ref{GUBk2}) and (\ref{RmUBk}), therefore
\begin{eqnarray*}
-\Delta|\nabla^{k-1}Rm|&\leq&C(n)(|Rm||\nabla^{k-1}Rm|+|\nabla
G||\nabla^{k+2}G|+|\nabla^{2}G||\nabla^{k+1}G|\\&&+C(\{\int|Rm|^{\frac{n}{2}}\}^{\frac{2}{n}}+\{\int|\nabla
  G|^2\}^{\frac{1}{2}})
\end{eqnarray*}
we have $L^{\frac{p}{2}} (p>n)$ bound of
$|\nabla^{2}G||\nabla^{k+1}G|+C(\{\int|Rm|^{\frac{n}{2}}\}^{\frac{2}{n}}+\{\int|\nabla
  G|^2\}^{\frac{1}{2}})$,
and $L^p$ bound of $\nabla^{k}Rm$ by (\ref{GWk2p}) and
(\ref{RmWkp}).  Apply lemma \ref{MIL} again, we obtain the
$L^\infty$ estimation of $\nabla^{k}Rm$.
 \qed

\section{Compactness of Harmonic - Einstein Equation}
In this section, we will give a sketch proof on the theorem
\ref{Compactness}, since the argument is very similar to the case of
Einstein metrics, \cite{Anderson89}, \cite{BKN89}, \cite{Tian90},
Bach flat metric with constant scalar curvature \cite{TV05b},
\cite{AAJV11}, K\"{a}hler Ricci soliton \cite{CS07}, and extremal
K\"{a}hler metric \cite{CW11}. By the way, the author had also
written a detailed proof for the removable singularity theorem in
the case of Bach flat metric with constant scalar curvature before
this work. As stated in the theorem, we have two aspects to show:
one is the convergence of Harmonic - Einstein equation in certain
topology, the other is smooth extension of the Harmonic - Einstein
equation across the singularity, which is called to be the removable
singularity theory.

First, with the assumption in  theorem \ref{Compactness}, we can
bound the energy and Sobolev constant, which is appeared in the
$\epsilon$ - regularity theorem \ref{eRegT}.

 \lemma With the assumption in theorem
\ref{Compactness}, there are constants $\Lambda_k, k=1,2,3$, which
are depending on $X, D, V, E$, but not on $i$, such that
\[\int_{M_i}|Rm(g_i)|^2\leq \Lambda_1, \; \int_{M_i}|dG_i|^2\leq \Lambda_2, \; C_S(M_i)\leq  \Lambda_3.\]

\proof In fact, C.Croke \cite{Croke80} proved that the isoperimetric
constant is bounded above by a constant depending only on a lower
bound for the Ricci curvature, lower bound on volume and an upper
bound on the diameter. In the later, based on Gromov's technique,
Anderson \cite{Anderson92} give a local version, which require on
local (Euclidean) volume growth condition. On the other hand,
isoperimetric constant is equal to Sobolev constant by Federer -
Fleming's theory.  In our case, $Ric=\lambda g+\langle dG,
dG\rangle\geq \lambda g$, $Diam \leq D$ and $Vol\geq V$, so we have
a uniform upper bound for the Sobolev constant: $C_S\leq C(D, V)$.

    With Sobolev constant, from lemma \ref{UBHF},  we have  $\sup_{M_i}|\nabla
G_i|\leq C(D, V, E)$. Moreover, Ricci curvature is two sided
bounded: $|Ric(g_i)|\leq C(D, V, E)$,  and the scalar curvature
$R=\lambda n+ |dG|^2$ is bounded too.  Recall the Gauss - Bonnet
formula on compact four manifold $M$, see \cite{Besse87},
\[\mathcal {X}(M)=\frac{1}{8\pi^2}\int_M|Rm|^2-|Ric-\frac{1}{4}Rg|^2\]
 then we have $\int_{M_i}|Rm|^2\leq C(X, D, V, E)$. \qed

Now, let us  give a sketch proof on the theorem \ref{Compactness}.

\proof Step I: As in the case of Einstein metric or Bach flat metric
with constant curvature, since we have established the local
regularity of Harmonic - Einstein equation, then the sequence will
converge as stated in the theorem by applying the Cheeger - Gromov
convergence, not only the convergence of the metric $g_i$, but also
with function $G_i$. More precisely, taking
$\epsilon=\epsilon(n,C_S)$ in theorem \ref{eRegT}, consider the sets
\[\mathcal {R}_i(r)=\{x\in M_i|\{\int_{B(x,r)}|Rm|^{\frac{n}{2}}\}^{\frac{n}{2}}<\epsilon\},\quad\mathcal {S}_i(r)=\{x\in M_i|\{\int_{B(x,r)}|Rm|^{\frac{n}{2}}\}^{\frac{n}{2}}\geq \epsilon\}\]
then $M_i=\mathcal {R}_i(r)\cup\mathcal {S}_i(r)$ and also $\mathcal
{R}_i(r_1)\subset\mathcal {R}_i(r_2)$, $\mathcal
{S}_i(r_1)\supset\mathcal {S}_i(r_2)$, for any $r_1 > r_2$.

For all $x\in\mathcal {R}_i(r)$, by $\epsilon$ - regularity theorem
\ref{eRegT},  for all $k\in\mathbb{N}$, we have
\[\sup_{B(0,\frac{r}{2})}|\nabla^{k} G_i|\leq\frac{C}{r^{k}},\quad
  \sup_{B(0,\frac{r}{2})}|\nabla^kRm(g_i)|\leq\frac{C}{r^{k+2}},\]
where $C=C(n,k,\Lambda_1,\Lambda_2,\Lambda_3)$.

Letting $\{B(x_k^i, \frac{r}{4})\}, k\in\mathbb{N}$ be a collection
of a maximal family of disjoint geodesic balls in $M_i$, then
$M_i\subset \cup_kB(x_k^i, r)$. There is a  uniform bound,
independent of $i$,  on the number of points $\{x_k^i\in
{S}_i(r)\}$, which follows from
\begin{equation}\label{FpSS}
m\leq\sum_{i=1}^m\epsilon^{-\frac{n}{2}}\int_{B(x_k^i,r)}|Rm|^{\frac{n}{2}}\leq
C\epsilon^{-\frac{n}{2}}\int_{M_i}|Rm|^{\frac{n}{2}},
\end{equation}
where $C=\sup_{x\in
M_i}\frac{Vol(B(x,\frac{5r}{4})}{Vol(B(x,\frac{r}{4})}\leq
C(n,\Lambda_2)$. Without loss generality, we will assume $m$ is
fixed, which is independent on $i$ and $r$.

On the other hand, the uniform Sobolev constant implies uniform
noncollapsing, namely, $Vol(B(x,r))\geq C(C_S)r^n$. Combine the
uniform bound of curvature, we have a uniform lower bound on the
local injective radius, i.e. $\hbox{inj}(x)\geq Cr$, $x\in
{R}_i(r)$, see \cite{CGT82}. According Cheeger - Gromov convergence
theory \cite{GW88}, we can extract a subsequence, so that $(\mathcal
{R}_j(r), g_j, G_j)$ converges smoothly to a smooth open Riemanniann
manifold $(\mathcal {R}_\infty(r),g_\infty, G_\infty)$. Since the
convergence is in the $C^\infty(\mathcal {R}_\infty(r))$ topology,
then the limit $(g_\infty, G_\infty)$ still satisfies the Harmonic -
Einstein equation on ${R}_\infty(r)$.

We now choosing a sequence $\{r_k\}\rightarrow 0$ and repeat the
above construction by choosing subsequence, we still denote $\{j\}$.
Since $\mathcal {R}_i(r_{k})\subset\mathcal {R}_i(r_{k+1})$, then we
have a sequence of limit spaces with natural inclusions
\[\mathcal{R}_\infty(r_{k})\subset\mathcal {R}_\infty(r_{k+1})\subset\cdots\subset\mathcal {R}_\infty:=dir.\lim\mathcal{R}_\infty(r_{k})\]

Due to finite capacity of $\mathcal {S}_j$ in (\ref{FpSS}),
following the argument of \cite{Anderson89}, \cite{Tian90}, one can
add finite points $\mathcal {S}_\infty=\{b_1,\cdots,b_m\}$ to
$\mathcal {R}_\infty$ such that $M_\infty:=\mathcal
{R}_\infty\cup\mathcal {S}_\infty$ is complete with respect
$g_\infty$. Since $|\nabla G|$ is uniformly bounded,
$G_\infty(b_k):=\lim_{b\rightarrow b_k}G_\infty(b)\in (N,h)$ is well
defined for $k=1,\cdots m$.

Moreover, with the local regularity, the curvature may blow up at
the singularity, but at worst, at a rate of quadratic, i.e.
$\sup_{\{x: d(x,\mathcal {S})=r\}}|Rm|\leq\frac{o(r)}{r^2}$, then we
know the singularity has a $C^0$ orbiflod structure, see
\cite{Tian90} or \cite{TV05a}. Since the energy is concentrated at
the singular set $\mathcal {S}$, then both $\int_{M_\infty}|Rm|^2$
and $\int_{M_\infty}|dG|^2$ remains bounded, which follows from the
lower semi - continuously of energy.

Step II: We will extend the Harmonic - Einstein equation across the
singularity. On the limit space with finite singularity, we may
assume the finite group is trivial, i.e. $\Gamma_k$=\{e\}, by going
to the universal covering space. If the full curvature is uniformly
bounded near the singularity, one can construct ``good" coordinate,
namely, $C^{1,\alpha}$ harmonic coordinate around the singularity
\cite{KD81}, \cite{BKN89}, \cite{Tian90}. And consequently, one can
go back to the equation and bootstrapping to improve the regularity.
Thus the main task is to bound the full curvature tensor.

We want to get the similar estimation via the Moser iteration on the
Riemannian orbifold. However, due to singularity of the manifold
structure, this would appear impossible, as we do not know that our
elliptic inequality $\Delta u\geq - f u- h$ holds weakly across the
singularity. We will easily see the Sobolev inequality does hold
despite the singularity, while integration by parts leaves an
uncontrollable term $\int|\nabla \phi|^2u^p$ near the singular
point. If $u\in L^{p}, p> \gamma $, let $\phi$ be zero on $B(p,r)$
near the singularity, then
\[\int|\nabla \phi|^2u^p\leq \{\int|\nabla \phi|^n\}^{\frac{n}{2}}\{\int_{supp|\nabla \phi|} |u|^{p\gamma}\}^{\frac{1}{\gamma}}\]
become negligible since
$\{\int_{B(p,r)}|\nabla\phi|^n\}^{\frac{2}{n}}$ is uniformally
bounded, this is Siber's lemma \cite{Sibner85}, which is used by
\cite{CS07} and \cite{CW11}.

When $n=4$, then $\gamma=2$; but we only have $u=|Rm|\in L^2$, thus
the equation does not hold in the weak sense across the singularity.
One approach to overcome this problem is using the Yang - Mill like
argument under Hodge gauge to improve the estimation, which is
created by Uhlenbeck \cite{Uhlenbeck82a}, and developed by Tian
\cite{Tian90}, \cite{TV05a}. More precisely, by choosing Hodge
gauge, integration on the annuls around the singularity point, one
can compare the energy of $Rm$, $f(r):=\int_{B(p_k,r)}|Rm|^2$, with
its derivatives, $f'(r)=\int_{S(p_k,r)}|Rm|^2$ to get a differential
inequality on $f(r)$. And consequently, one can improve the decay
order of $f(r)$.

\lemma\cite{Tian90}\label{Hodge gauge} There are constants
$\epsilon$ and $C$ such that any connection $A$ on the trivial
bundle over a punctured ball $B(0,1)\setminus \{0\}$ with
$\|R_A\|\leq \frac{\epsilon(r)}{r^2}$, is gauge equivalent to a
connection $A^\tau$ on the annulus $\Omega(r,R):=B(0,R)\setminus
\overline{B(0,r)}$ with
\begin{enumerate}
  \item $d^*A^\tau(r,R)=0 \quad \hbox{in}\quad \overline{\Omega(r,R)}$,
  \item $d_{\psi}^*A^\tau=0\quad \hbox{on}\quad S(r,R):=\partial \Omega(r,R)$,
  \item $\int_{\Omega(r,R)}A(\nabla r)=0$,
  \item $\int_{\Omega(r, 2r)}|A|^2\leq Cr^2\int_{\Omega(r, 2r)}|R_A|^2$
\end{enumerate}
where $d^*$ and $d_{\psi}^*$ are the adjoint operators of the
exterior differentials on $\Omega(r,R)$, $S(r,R)$ respectively.
Moreover, for suitable constants $\epsilon$ and $C$, the connection
$A^\tau$ is uniquely determined, up to the transformation
$A^\tau\rightarrow u_0A^\tau u_0^{-1}$ for constant gauge
$u_0$.\footnote{It follows from the uniqueness of Hodge gauge on
sphere, Theorem 2.5 \cite{Uhlenbeck82a}.}

Now we will improve the decay order of $\int |Rm|^2$ by the same
argument in \cite{Tian90} or \cite{TV05a}, but change only a few
words, namely, the Ricci term is related to the harmonic map $G$,
see (\ref{W2L2aS}).

\lemma \label{CDaS} For any $b_k\in\mathcal {S}_\infty$, denote
$B(r):=B(b_k,r)$ and $S(r):=\partial B(b_k,r)$. There exists
$1<\beta<2$ such that for $r$ sufficiently small, we have
\begin{equation*}
    \sup_{S(r)}|Rm|\leq C r^{-(2-\beta)}.
\end{equation*}
\proof Choose $r_0=r$ small, and let us denote
$r_i=\frac{1}{2}r_{i-1}$. Let $A_i$ be the connection on $\Omega_i
=\Omega(r_i, r_{i-1})$ from lemma \ref{Hodge gauge}, then there
exist Hodge Gauge
\begin{eqnarray*}
  d_{\psi}^*A_{i\psi}\mid_{\partial \Omega_i} &=& 0 \\
  d_{\psi}^*A_{(i+1)\psi}\mid_{\partial \Omega_{i+1}} &=& 0
\end{eqnarray*}
so the restriction $A_{i\psi}$ and $A_{(i+1)\psi}$ differ by a
constant gauge on $S(r_i)$ and we may therefore assume
that\[A_{i\psi}\mid_{S(r_i)}=A_{(i+1)\psi}\mid_{S(r_i)}\] and then
the curvature is continuous across the $S_i$, i.e.
$(R_{A_i})_{r\psi}=(R_{A_{i+1}})_{r\psi}$ follows from the gauge
transformation rule of curvature. Then we compute the $L^2$ of
curvature
\begin{eqnarray*}
  \int_{\Omega_i}|Rm|^2 &=&\int_{\Omega_i}\langle D_i A_i - [A_i, A_i], R_{A_i}\rangle \\
   &=& -\int_{\Omega_i} \langle A_i, D_i^*R_{A_i}\rangle-\int_{\Omega_i}\langle [A_i, A_i], R_{A_i}\rangle \\
&&- \int_{S_i}\langle  (A_i)_\psi,
(R_{A_i})_{r\psi}\rangle-\int_{S_{i+1}}\langle  (A_i)_\psi,
(R_{A_i})_{r\psi}\rangle
 \end{eqnarray*}
Next we sum over $i$, the boundary terms cancel, except for $S_0$
and the inner budgetary terms become negligible as $i\rightarrow
\infty$,
\begin{eqnarray*}
  \int_{B(r)}|Rm|^2 
   &=& \sum_{i=1}^\infty\int_{\Omega_i}|R_{A_i}|^2\\
   &=& -\sum_{i=1}^\infty\int_{\Omega_i} \langle A_i, D_i^*R_{A_i}\rangle-\sum_{i=1}^\infty\int_{\Omega_i}\langle [A_i, A_i], R_{A_i}\rangle \\
&&+ \int_{S(r)}\langle  (A_1)_\psi, (R_{A_1})_{r\psi}\rangle
 \end{eqnarray*}
Let us estimate the three term on the right separatively. In fact,
on the round sphere $S^3$, the first eigenvalue for the Laplacian on
coclosed 1 form is 4 \footnote{See also the Corollary 2.6 of
\cite{Uhlenbeck82a}.}, then
\begin{eqnarray*}
  \{4\int_{S^3}|A|^2\}^{\frac{1}{2}} &\leq& \{\int_{S^3}|dA|^2\}^{\frac{1}{2}} \\
  &\leq& \{\int_{S^3}|R_A|^2\}^{\frac{1}{2}}+C|R_A|_{L^\infty}\{\int_{S^3}|A|^2\}^{\frac{1}{2}}
\end{eqnarray*}
Since the singularity is $C^0$ orbifold, then the geodesic sphere is
convergence to the round sphere after scaling. Therefore we may find
monotone function $\epsilon'(r)$ with $\lim_{r\rightarrow
0}\epsilon'(r)=0$ such that
\begin{eqnarray*}
\int_{S(r)}\langle  (A_1)_\psi, (R_{A_1})_{r\psi}\rangle &\leq& (\int_{S(r)}|(A_1)_\psi|^2)^{\frac{1}{2}}(\int_{S(r)}|(R_{A_1})_{r\psi}|^2)^{\frac{1}{2}}\\
 &\leq&\frac{1}{2-\epsilon'(r)} r(\int_{S(r)}|(R_{A_1})_{\psi\psi}|^2)^{\frac{1}{2}} (\int_{S(r)}|(R_{A_1})_{r\psi}|^2)^{\frac{1}{2}}\\
&\leq&\frac{1}{2}\frac{1}{2-\epsilon'(r)}r\int_{S(r)}|Rm|^2
\end{eqnarray*}
Note that $D_i^*R_{A_i}= \nabla^*Rm=d^{\nabla}Ric$,  and we estimate
the first term as
\begin{eqnarray*}
  \int_{\Omega_i} \langle A_i, D_i^*R_{A_i}\rangle &\leq& \frac{\delta}{C} r_i^{-2}\int_{\Omega_i}|A_i|^2+ C\delta^{-1}r_i^2\int_{\Omega_i}|\nabla Ric|^2\\
  &\leq& \delta\int_{\Omega_i}|R_{A_1}|^2+C\delta^{-1}r_i^2\int_{\Omega_i}|\nabla Ric|^2
\end{eqnarray*}
where the $C$ is the constant in $4$ - th item of lemma \ref{Hodge
gauge}. Moreover, we have
\begin{eqnarray*}
  |\sum_{i=1}^\infty\int_{\Omega_i}\langle [A_i, A_i], R_{A_i}\rangle| &\leq&\sum_{i=1}^\infty C\sup_{\Omega_i}|R_{A_i}|\int_{\Omega_i}|A_i|^2 \\
   &\leq& \sum_{i=1}^\infty C\epsilon(r_{i-1})\int_{\Omega_i}|R_{A_i}|^2\\
   &=&C\epsilon(r)\int_{B(p,r)}|Rm|^2
\end{eqnarray*}
For the $Ricci$ term,  we will prove later that we still have
(\ref{G2L2}),
\begin{equation}\label{W2L2aS}
\int_{\Omega_i}|\nabla Ric|^2\leq Cr_i^2
\end{equation}
Combining the above stimulation, we obtain
\[(1-C\epsilon(r)-\delta)\int_{B(r_0)}|Rm|^2\leq \frac{1}{2}\frac{1}{2-\epsilon'(r)}\int_{S(r)}|Rm|^2+ \sum_{i=1}^\infty C\delta^{-1}r_i^{4}\]
Therefore for all $r$ sufficiently small, choosing $\delta$
sufficiently small, there exists a small constant $\delta'\geq 0$,
\[\int_{B(p,r)}|Rm|^2\leq \frac{r}{4-\delta'}\int_{S(r)}|Rm|^2+Cr^{4}\]
We denote $f(r)=\int_{B(p,r)}|Rm|^2$, then
\[f(r)\leq \frac{1}{c_1}rf'(r)+c_2r^{4}\]
where  $c_1=4-\delta'$. Then
\[(r^{-c_1}f(r))'=r^{-c_1}f'(r)-c_1r^{-c_1-1}f(r)\geq-c_1c_2r^{-1+\delta'},\]
Therefore
\[r_0^{-c_1}f(r_0)-r^{-c_1}f(r)=\int_r^{r_0}(r^{-c_1}f(r))'\geq-c_1c_2\int_r^{r_0}r^{-1+\delta'}\]
If $\delta'> 0$, then $f(r)\leq
Cr^{4-\delta'}+\frac{C}{\delta'}r^4\leq Cr^{4-\delta'}$; If
$\delta'=0$, then $f(r)\leq Cr^{4}+Cr^4\ln r$. In any case,
$\delta'$ is small, so we can choose $1<\beta<2$ such that $f(r)\leq
Cr^{2\beta}$. Therefore, working on the ball
$B(x,\frac{r}{2})\subset B(p,1)\setminus\{p\}$, $x\in S(r)$, from
the $\epsilon$- regularity theorem \ref{eRegT} for the smooth case,
we have
\[|Rm|(x)\leq\sup_{B(x,\frac{r}{4})}|Rm| \leq C(n) r^{-2}\{C_S^n\int_{B(x,\frac{r}{2})} |Rm|^2\}^\frac{1}{2}\leq Cr^{-2+\beta}.\]\qed

\proof (proof of (\ref{W2L2aS})) Define $\phi_l=\eta_l\phi$, where
$\eta_l\equiv 0$ in $B(b_k,\frac{1}{2l})$, $\eta_l\equiv 1$ in
$B(b_k,r)\setminus B(b_k,\frac{1}{l})$, $|\nabla\eta_l|\leq 4l$ when
$l$ is large. Recall the proof of (\ref{G2L2}), let
$\phi_l=\eta_l\phi$ replace $\phi$ as to be the cut off function on
$B(b_k, r)$, since $\phi_l$ vanish at the singular point $b_k$, as
in the smooth case, we have (\ref{G2L2}):
\[\int(\phi_l)^2|\nabla^2G|^2\leq C(\int (\phi_l)^2|\nabla G|^2+\int |\nabla(\phi_l)|^2|\nabla G|^2)\]
On the other hand,
\begin{eqnarray*}
\int(\phi_l)^2|\nabla^2G|^2&\leq& C(\int (\phi_l)^2|\nabla G|^2+\int |\nabla(\phi_l)|^2|\nabla G|^2) \\
&\leq & C(\int \phi^2|\nabla G|^2+\int |\nabla\phi|^2|\nabla G|^2+\int|\nabla\eta_l|^2|\nabla G|^2 ) \\
&\leq & C(\int \phi^2|\nabla G|^2+\int |\nabla\phi|^2|\nabla
G|^2+\sup|\nabla G|^2\frac{1}{l^2} )
\end{eqnarray*}
Letting $l$ tends to $\infty$, so we have
\begin{eqnarray*}
\int\phi^2|\nabla^2G|^2&\leq & C(\int \phi^2|\nabla G|^2+\int
|\nabla\phi|^2|\nabla G|^2)
\end{eqnarray*}
and
\[\int_{B(b_k,\frac{r}{2})}|\nabla^2G|^2\leq
C\frac{1}{r^2}\int_{B(b_k,r)}|\nabla G|^2\leq Cr^2.\]\qed

With lemma \ref{CDaS}, even though we do not have the uniform
bounded curvature across the singularity, the curvature condition
$\sup_{S(r)}|Rm|\leq C r^{-(2-\beta)}$ with $1<\beta<2$ is enough to
construct $C^{1,\beta-1}$ coordinate around the singularity, see
\cite{BKN89} and \cite{Tian90}. More precisely, we can construct
coordinates $\Psi: S^3\times (0, 1]\rightarrow \widetilde{B}(b_k,
1)$ such that,
\[\Psi^*g_{ij}(x)-\delta_{ij}=O(|x|^\beta),\quad
\partial\Psi^*g_{ij}(x)=O(|x|^{\beta-1}).\]

By \cite{KD81}, one can  construct harmonic coordinates around the
singularity with regularity at least $C^{1,\alpha}$,
$\alpha=\beta-1$. Now apply the Schauder theory on the coupled
system (\ref{HaMEl})
\begin{equation}
\left\{ \begin{aligned}
\Delta g&=-2Ric+Q(g,\partial g)=-2\langle dG, d G\rangle +Q(g,\partial g)\\
\Delta G&=dG\ast dG
\end{aligned} \right.
 \end{equation}
under the harmonic coordinate.  We first have that $g\in
C^{1,\alpha}$ and $dG\in \L^\infty$. By the second equation, $G\in
C^{1,\alpha}$; going back to the first equation, the right hand side
is $C^{\alpha}$, therefore $g\in C^{2,\alpha}$; Go to the second
equation again, $G\in C^{2,\alpha}$; and consequently, by the first
equation again, $g\in C^{3,\alpha}$. Bootstrapping in this manner,
we actually show that both $g$ and $G$ is smooth across the
singularity.  So we finish the proof of theorem
\ref{Compactness}.\qed

\end{document}